\documentclass[11pt]{paper}

\usepackage{cite}
\usepackage{pifont}
\usepackage{times}
\usepackage{amssymb,epsfig,floatflt}
\usepackage{subfigure}[1995/03/06]
\usepackage[section]{placeins}
\usepackage{epsfig,psfrag}
\usepackage{graphics}
\usepackage{color}

\input amssymb.sty

\newcommand{\IR}{\mathbb{R}}
\newcommand{\II}{\mathbb{I}}

\newcommand{\sign}{\mathrm{sgn}}

\newtheorem{thm}{Theorem}

\newtheorem{lem}{Lemma}

\newtheorem{defn}{Definition}

\newtheorem{cor}{Corollary}

\newtheorem{rem}{Remark}

\usepackage{cite}

\begin{document}
\thispagestyle{empty}

\title{\LARGE \bf Design and Analysis of a Novel $\mathcal{L}_1$ Adaptive Control Architecture with Guaranteed Transient Performance\thanks{Research is supported by AFOSR under Contract
No. FA9550-05-1-0157.}}
\author{Chengyu Cao and Naira Hovakimyan\thanks{The authors are with
  Aerospace \& Ocean Engineering, Virginia Polytechnic
 Institute \& State University, Blacksburg, VA 24061-0203, e-mail:
  chengyu, nhovakim@vt.edu}}
\date{}
\maketitle
\begin{abstract}

 Conventional Model Reference Adaptive Controller (MRAC),
while providing an architecture for  control of systems in the
presence of parametric uncertainties, offers no means for
characterizing the system's input/output performance during the
transient phase. Application of adaptive controllers was therefore
largely restricted due to the fact that the system uncertainties
during the transient have led to unpredictable/undesirebale
situations, involving control signals of high-frequency or large
amplitudes, large transient errors or slow convergence rate of
tracking errors, to name a few.
 In this paper, we develop a novel  adaptive
control architecture that ensures that the  input and the output
of an uncertain linear system  track the input and output of a
desired linear system during the transient phase, in addition to
the asymptotic tracking. This new architecture has a low-pass
filter in the feedback-loop that enables to enforce the desired
transient performance by increasing the adaptation gain. For the
proof of asymptotic stability, the $\mathcal L_1$ gain of a
cascaded system, comprised of this filter and the closed-loop
desired reference model, is required to be less than the inverse
of the upper bound of the norm of unknown parameters used in
projection based adaptation laws.
The ideal (non-adaptive)
version of this $\mathcal L_1$ adaptive controller
   is used along with the main system dynamics to define a closed-loop reference system, which gives an
opportunity to estimate performance bounds in terms of $\mathcal
L_{\infty}$ norms  for both  system's input and output signals as
compared to the same signals of this reference system.  Design
guidelines  for selection of the low-pass filter  ensure that the
 closed-loop reference system  approximates the desired
system response, despite the fact that it depends upon the unknown
parameters.  The tools from this paper can be used to develop a
theoretically justified verification and validation framework for
adaptive systems. Simulation results illustrate the theoretical
findings.
\end{abstract}

\section{Introduction}

Model Reference Adaptive Control (MRAC) architecture was developed
conventionally to control  linear systems in the presence of
parametric uncertainties \cite{NarBook89,SloBook91}. The
development of this architecture has been facilitated by the
Lyapunov stability theory that defines sufficient conditions for
stable performance, but offers no means for characterizing the
 system's
input/output performance during the transient phase. Application
of adaptive controllers was therefore largely restricted due to
the fact that the system uncertainties during the transient have
led to unpredictable/undesirebale situations, involving control
signals of high-frequency or large amplitudes, large transient
errors or slow convergence rate of tracking errors, to name a few.
Extensive tuning of adaptive gains and Monte-Carlo runs have been
the primary methods up today enabling the transition of  adaptive
control solutions to real world applications. This argument has
rendered verification and validation of adaptive controllers
overly challenging.

Improvement of the transient  performance of adaptive controllers
has been addressed from various perspectives in numerous
publications \cite{SloBook91, DattaHo94, BartFerSto99, Sun93,
MilTac91acw, Costa99, IoaBook03, YdsTac92tpa, KrsScl93tpi,
OrtTac93mna, ZanCdc90tbf, Datta94, arteaga2002, narendra94}, to
name a few. An example  presented in \cite{ZanCdc90tbf}
demonstrated that the system output can have overly poor transient
tracking behavior before ideal asymptotic convergence can take
place. On the other hand, in \cite{YdsTac92tpa} the author proved
that it may not be possible to optimize $\mathcal{L}_2$ and
$\mathcal{L}_{\infty}$ performance simultaneously by using a
constant adaptation rate. Following these results,  modifications
of adaptive controllers were proposed in \cite{Sun93, Datta94}
that render the tracking error arbitrarily small in terms of both
mean-square and $\mathcal{L}_{\infty}$ bounds.  Further, it was
 shown in \cite{DattaHo94} that the
 modifications proposed in \cite{Sun93, Datta94}  could be
derived as a linear feedback of the tracking error, and the
 improved performance was obtained
 only due to a nonadaptive high-gain feedback in that
scheme.
 In \cite{SloBook91}, composite
adaptive controller was proposed, which suggests a new adaptation
law using both tracking error and prediction error that leads to
less oscillatory behavior in the presence of high adaptation gains
as compared to MRAC.  In \cite{MilTac91acw},  a high-gain
switching MRAC technique was introduced to achieve arbitrary good
transient tracking performance under a relaxed set of assumptions
as compared to MRAC, and the results were shown to be of existence
type only.  In \cite{narendra94},  multiple model and switching
scheme is proposed to improve the transient performance of
adaptive controllers. In \cite{arteaga2002}, it is shown that
arbitrarily close transient bound can be achieved by enforcing
parameter-dependent persistent excitation condition.
 In \cite{KrsScl93tpi},
computable $\mathcal{L}_2$ and $\mathcal{L}_{\infty}$ bounds for
the output tracking error signals are obtained for a special class
of adaptive controllers using backstepping. The underlying linear
nonadaptive controller possesses a parametric robustness property,
however, for a large parametric uncertainty it requires high gain.
In \cite{OrtTac93mna}, dynamic certainty equivalent controllers
with unnormalized estimators were used for adaptation that permit
to derive a uniform upper bound for the $\mathcal{L}_2$ norm of
the tracking error in terms of initial parameter estimation error.
In the presence of sufficiently small initial conditions, the
author proved that the $\mathcal{L}_{\infty}$ norm of the tracking
error is upper bounded by the $\mathcal{L}_{\infty}$ norm of the
reference input. In \cite{PanBasar98, ArslanBasar2001},
differential game theoretic approach has been investigated for
achieving arbitrarily close disturbance attenuation for tracking
performance, albeit at the price of increased control effort. In
\cite{ioannou2000}, a new certainty equivalence based adaptive
controller is presented using backstepping based control law with
a normalized adaptive law to achieve asymptotic stability and
guarantee performance bounds comparable with the tuning functions
scheme, without the use of higher order nonlinearities.

As compared to the linear systems theory, several important
aspects of the transient performance analysis seem to be missing
in these papers. First, all the bounds  in these papers are
computed for tracking errors only, and not for  control signals.
Although the latter can be deduced from the former, it is
straightforward to verify that the ability to adjust the former
may not extend to the latter in case of nonlinear control laws.
Second, since the purpose of adaptive control is to ensure stable
performance in the presence of modeling uncertainties, one needs
to ensure that the changes in reference input and unknown
parameters due to possible faults or unexpected uncertainties do
not lead to unacceptable transient deviations or oscillatory
control signals, implying that a retuning of adaptive parameters
is required. Finally, one needs to ensure that whatever
modifications or solutions are suggested for performance
improvement of adaptive controllers, they are not achieved via
high-gain feedback.

Following \cite{kkk}, subject to appropriate trajectory
initialization, the following bound $ ||e||_{\infty}\le V(t)\le
\frac{V(0)}{\lambda_{\min}(P)}\le\frac{\tilde \theta^2(0)}{\Gamma}
$, where $e$ is the tracking error,  $\tilde \theta$ is the
parametric error, $V(t)$ is the positive definite Lyapunov
function, $\lambda_{\min}(P)$ is the minimum eigenvalue of
$P=P^\top>0$, found by solving the algebraic Lyapunov equation
associated with the error dynamics, implies that increasing the
adaptation gain $\Gamma$ leads to smaller tracking error for all
$t\ge 0$, including the transient phase. However, large adaptive
gain leads to high frequencies in the control signal, implying
that the improvement in the transient tracking of the system
output is achieved   at the price of unacceptable high frequencies
in the system input. One can observe from the open-loop transfer
function analysis for a PI controller (which is a MRAC-structure
controller for a linear system with constant disturbance) that
increasing the adaptation gain leads to reduced phase-margin, and
consequently reduced time-delay tolerance in input/output
channels. On the contrary, decreasing the adaptive gain leads to
large deviations from the desired trajectory during the transient
phase.

 In this paper we define a new type of model following adaptive
 controller that
  adapts fast leading to
 desired transient performance  for  the
 system's  both input and output signals simultaneously.  This new architecture has a low-pass filter
in the feedback-loop that enables to enforce the desired transient
performance by increasing the adaptation gain. For the proof of
asymptotic stability, the $\mathcal L_1$ gain of a cascaded
system, comprised of this filter and the closed-loop desired
transfer function, is required to be less than the inverse of the
upper bound on the norm of unknown parameters used in projection
based adaptation laws. With the  low-pass filter in the loop, the
$\mathcal L_1$ adaptive controller is guaranteed to stay in the
low-frequency range
  even in the presence of
   high adaptive gains and large reference inputs.  The ideal (non-adaptive) version of this $\mathcal L_1$ adaptive controller
   is used along with the main system dynamics to define a  closed-loop reference system, which gives an
opportunity to estimate performance bounds in terms of $\mathcal
L_{\infty}$ norms  for  system's both input and output signals as
compared to the same signals of this reference system. These
bounds immediately imply that the transient performance of the
control signal in MRAC cannot be characterized. Design guidelines
for selection of the low-pass filter  ensure that the closed-loop
reference system  approximates the desired system response,
despite the fact that it depends upon the unknown parameter. Thus,
the desired tracking performance is achieved by systematic
selection of the low-pass filter, which in its turn enables  fast
adaptation, as opposed to high-gain designs leading to increased
control efforts. Using a simple linear system with constant
disturbance, we demonstrate that this new architecture has
guaranteed time-delay margin in the presence of large adaptive
gain, as opposed to MRAC. We further demonstrate the extension of
the methodology to systems with unknown time-varying parameters.

The paper is organized as follows. Section \ref{sec:preliminary}
states some preliminary definitions, and Section \ref{sec:PF}
gives the problem formulation.  In Section \ref{sec:Hinf}, the new
$\mathcal{L}_1$ adaptive controller is presented. Stability and
tracking results of the $\mathcal{L}_1$ adaptive controller are
presented in Section \ref{sec:convergence}. Design guidlines are
provided in Section \ref{sec:design}. Comparison of the
performance of $\mathcal{L}_1$ adaptive controller, MRAC and the
high gain controller are discussed in section \ref{sec:comp}.
Analysis of $\mathcal{L}_1$ adaptive controller in  the presence
of time-varying unknown parameters is presented in Section
\ref{sec:timevarying}. In Section \ref{sec:simu}, simulation
results are presented, while Section \ref{sec:con} concludes the
paper. Unless otherwise mentioned, the $||\cdot||$ will be used
for the $2$-norm of the vector.

\section{Preliminaries}\label{sec:preliminary}

In this Section, we recall  basic definitions and facts from
linear systems theory, \cite{IoaBook03,KhaBook02,ZhoBook98}.

%
%
%
%
%
%
%
%
%
%

\begin{defn}\label{defn1}
For a signal $\xi(t), t\geq 0, \xi\in {\IR}^n$, its truncated
${\mathcal L}_{\infty}$ norm and ${\mathcal L}_{\infty}$ norm are
$\displaystyle{ \Vert \xi_t \Vert_{{\mathcal L}_\infty} =
\max_{i=1,..,n} (\sup_{0\leq \tau \leq t} |\xi_i(\tau)|)}$,
$\displaystyle{\Vert \xi \Vert_{{\mathcal L}_\infty} =
\max_{i=1,..,n} (\sup_{\tau\geq 0} |\xi_i(\tau)|),}$ where $\xi_i$
is the $i^{th}$ component of $\xi$.
\end{defn}
\begin{defn}
 The $\mathcal{L}_1$ gain of a stable proper single--input single--output system $H(s)$
 is defined to be
$|| H(s)||_{\mathcal{L}_1} = \int_{0}^{\infty} |h(t)| d t$.
\end{defn}

{\em Proposition:}
 A continuous time LTI
system (proper) with impulse response $h(t)$ is stable if and only
if
$
\int_{0}^{\infty} |h(\tau)| d\tau < \infty.
$
A proof can be found in \cite{IoaBook03} (page 81, Theorem 3.3.2).

\begin{defn}
For a stable proper $m$ input $n$ output system $H(s)$ its
$\mathcal{L}_1$ gain is defined as
\begin{equation}\label{L1def}
\Vert H(s) \Vert_{\mathcal{L}_1} = \max_{i=1,..,n} (
\sum_{j=1}^{m}\Vert H_{ij}(s)\Vert_{\mathcal{L}_1} )\,,
\end{equation}
where $H_{ij}(s)$ is the $i^{th}$ row $j^{th}$ column element of
$H(s)$.
\end{defn}

The next lemma extends the results of Example 5.2.
(\cite{KhaBook02}, page 199)  to general multiple input multiple
output systems.

\begin{lem}\label{lem:L1} For a stable proper multi-input multi-output (MIMO)  system $H(s)$ with input $r(t) \in {\IR}^m$ and output
 $x(t)\in {\IR}^n$, we have
\begin{equation}
\Vert x_t\Vert_{{\mathcal L}_\infty} \leq \Vert H\Vert_{{\mathcal
L}_1} \Vert r_t\Vert_{{\mathcal L}_\infty},\quad \forall~
t>0.\nonumber
\end{equation}
\end{lem}

\begin{cor}\label{lem:L1_ext} For a stable proper MIMO  system  $H(s)$, if the input   $r(t) \in {\IR}^m$ is bounded, then the  output
 $x(t)\in {\IR}^n$ is also bounded as
$ 
\Vert x\Vert_{{\mathcal L}_\infty } \leq \Vert
H(s)\Vert_{\mathcal{L}_1} \Vert r\Vert_{{\mathcal L}_\infty }.
$ 
\end{cor}

\begin{lem}\label{lem:L1cas} For a cascaded system $H(s)=H_2(s) H_1(s)$,
where $H_1(s)$ is a stable proper system with $m$ inputs and $l$
outputs
 and $H_2(s)$ is a stable   proper system with $l$ inputs and  $n$
outputs, we have
$ 
\Vert H(s) \Vert_{\mathcal{L}_1} \leq  \Vert H_2(s)
\Vert_{\mathcal{L}_1} \Vert H_1(s) \Vert_{\mathcal{L}_1}\,.
$
\end{lem}
{\bf Proof.} Let $y(t)\in {\IR}^n$ be the  output of $H(s)=H_1(s)
H_2(s)$ in response to input $r(t)\in {\IR}^m$. It follows from
Lemma \ref{lem:L1_ext} that
\begin{equation}\label{lemL1cas6}
\Vert y(t) \Vert \leq ||y||_{{\mathcal L}_{\infty}} \leq \Vert
H_2(s) \Vert_{\mathcal{L}_1} \Vert H_1(s) \Vert_{\mathcal{L}_1}
\Vert r\Vert_{{\mathcal L}_\infty}
\end{equation}
for any bounded $r(t)$. Let $H_i(s), i=1,..,n$ be the $i^{th}$ row
of the system $H(s)$. It follows from (\ref{L1def}) that there
exists $i$ such that
\begin{equation}\label{lemL1cas9}
\Vert H(s) \Vert_{\mathcal{L}_1} = \Vert H_i(s)
\Vert_{\mathcal{L}_1}.
\end{equation}
Let $h_{ij}(t)$ be the $j^{th}$ element of the impulse response of
the system $H_i(s)$. For any $T$, let
\begin{equation}\label{lemL1cas7}
r_j(t) = \sign h_{ij}(T-t), \qquad t\in [0, T]\, , \;\forall
j=1,..,m.
\end{equation}
It follows from Definition \ref{defn1} that
$ 
\Vert r \Vert_{\mathcal{L}_{\infty}} =1\,,
$ 
and hence
$ 
\Vert y(t) \Vert \leq  \Vert H_2(s) \Vert_{\mathcal{L}_1} \Vert
H_1(s) \Vert_{\mathcal{L}_1}\,,\forall t\geq 0.
$ 
 For $r(t)$ satisfying
(\ref{lemL1cas7}), we have
$$ 
y(T) = \int_{t=0}^{T} \sum_{j=1}^{m} h_{ij}(T-t) r_j(t) d t =
\int_{t=0}^{T} \sum_{j=1}^{m} | h_{ij}(T-t)| d t = \sum_{j=1}^{m}
 ( \int_{t=0}^{T} |h_{ij}(t)| d t).
$$ 
Therefore, it follows from (\ref{lemL1cas6}) that for any $T$,
$ 
 \sum_{j=1}^{m}
 \left( \int_{t=0}^{T} |h_{ij}(t)| d t \right)\leq \Vert H_2(s) \Vert_{\mathcal{L}_1} \Vert H_1
 \Vert_{\mathcal{L}_1}\,.
$ 
 As $T\rightarrow \infty$, it follows from (\ref{lemL1cas9}) that
$$ 
\Vert H(s) \Vert_{\mathcal{L}_1} =\Vert H_i(s)
\Vert_{\mathcal{L}_1} = \lim_{T\rightarrow \infty} \sum_{j=1}^{m}
( \int_{t=0}^{T} |h_{ij}(t)| d t) \leq \Vert H_2(s)
\Vert_{\mathcal{L}_1} \Vert H_1(s) \Vert_{\mathcal{L}_1}\,,
$$ 
and this completes the proof. $\hfill{\square}$

\begin{figure}[!t]
\begin{center}
\includegraphics[width=2.5in,height=1.4in]{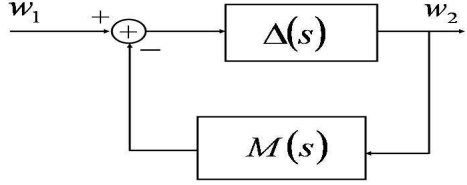}
\caption{Interconnected systems}\label{fig:sgblock}
\end{center}
\end{figure}

 Consider the interconnected
LTI system in Fig. \ref{fig:sgblock}, where $w_1\in {\IR}^{n_1}$,
$w_2\in {\IR}^{n_2}$, $M(s)$ is a stable  proper system with $n_2$
inputs and $n_1$ outputs, and $\Delta(s)$ is a stable proper
system with $n_1$ inputs and $n_2$ outputs.
\begin{thm}\label{thm:sg}
({\bf $\mathcal{L}_1$ Small Gain Theorem})  The interconnected
system in Fig. \ref{fig:sgblock} is stable if
$ 
\Vert M(s) \Vert_{\mathcal{L}_1} \Vert \Delta(s)
\Vert_{\mathcal{L}_1} < 1.
$ 
\end{thm}
The proof follows from Theorem 5.6  (\cite{KhaBook02}, page 218),
written for $\mathcal{L}_1$ gain.

%

Consider a linear time invariant system:
\begin{equation}\label{LTI_def}
\dot{x}(t) = A x(t) +b u(t)\,,
\end{equation}
where $x\in {\IR}^n$, $u\in {\IR}$, $b\in {\IR}^{n}$, $A\in
{\IR}^{n \times n}$ is Hurwitz, and assume that the transfer
function $(s I -A)^{-1} b $ is strictly proper and stable. Notice
that it can be expressed as:
\begin{equation}\label{LTI_trans}
(s I -A)^{-1} b = n(s)/d(s)\,,
\end{equation}
where $d(s)={\rm{det}} (s I- A)$ is a $n^{th}$ order stable
polynomial, and $n(s)$ is a $n\times 1$ vector with its $i^{th}$
element being a polynomial function:
\begin{equation}\label{nijdef}
n_i(s) = \sum_{j=1}^{n} n_{ij} s^{j-1} \,
\end{equation}

\begin{lem}\label{lem:pre1}
If $(A\in {\IR}^{n\times n}, b\in {\IR}^{n})$ is controllable, the
matrix $N$ with its $i^{th}$ row $j^{th}$ column entry $n_{ij}$ is
full rank.
\end{lem}
{\bf Proof.} Controllability of $(A, b)$ for the LTI system in
(\ref{LTI_def}) implies that given an initial condition $x(t_0)=0$
and  arbitrary $x_{t_1}\in {\IR}^n$ and arbitrary $t_1$, there
exists $u(\tau), \tau\in [t_0,\;t_1]$ such that $
x(t_1)=x_{t_1}\,. $ If $N$ is not full rank, then there exists a
non-zero vector $ u \in{\IR}^n$, such that $ u^{\top} n(s) =0\,. $
Then it follows that for $x(t_0)=0$ one has $ u^{\top}(\tau)
x(\tau) =0, ~ \forall \tau>t_0\,. $ This contradicts
$x(t_1)=x_{t_1}$, in which $x_{t_1}\in {\IR}^n$ is assumed to be
an arbitrary point. Therefore, $N$ must be full rank, and the
proof is complete. $\hfill{\square}$

\begin{lem}\label{lem:pre2}
If $(A, b)$ is controllable and $(s I -A)^{-1} b $ is strictly
proper and stable, there exists $c\in {\IR}^n$ such that the
transfer function $c^{\top}(s I -A)^{-1} b$ is minimum phase
 with relative degree one, i.e. all its zeros are located in
the left half plane, and its denominator is one order larger than
its numerator.
\end{lem}
{\bf Proof.} It follows from (\ref{LTI_trans}) that
\begin{equation}\label{lem_pre2_3}
c^{\top}(s I -A)^{-1} b = (c^{\top} N [s^{n-1}~ \cdots ~
1]^{\top})/d(s),
\end{equation}
where $N\in {\IR}^{n\times n}$ is matrix with its $i^{th}$ row
$j^{th}$ column entry $n_{ij}$ introduced in (\ref{nijdef}). We
choose $\bar{c}\in {\IR}^n$ such that $\bar{c}^{\top} [s^{n-1}~
\cdots ~ 1]^{\top}$ is a stable $n-1$ order polynomial. Since $(A,
b)$ is controllable, it follows from Lemma \ref{lem:pre1} that $N$
is full rank. Let $ c= (N^{-1})^{\top} \bar{c}. $ Then it follows
from (\ref{lem_pre2_3}) that $ c^{\top}(s I -A)^{-1} b =
\frac{\bar{c}^{\top} [s^{n-1}~ \cdots ~ 1]^{\top}}{d(s)} $ has
relative degree $1$ with all its zeros  in the left half plane.
$\hfill{\square}$

\section{Problem Formulation}\label{sec:PF}
Consider the following single-input single-output system dynamics:
\begin{eqnarray}
\dot{x}(t) & = & A x(t)+b (u(t)-\theta^{\top} x(t)),\quad x(0)=x_0 \label{problem}\\
y(t) & = & c^{\top} x(t)\,, \nonumber
\end{eqnarray}
 where $x\in
{\IR}^n$ is the system state vector (measurable), $u\in {\IR}$ is
the control signal, $b,c\in {\IR}^n$ are known constant vectors,
$A$ is known $n \times n$ matrix, $(A,b)$ is controllable, unknown
parameter $\theta\in \IR^n$ belongs to a given compact convex set
$ \theta\in \Omega$, $y\in {\IR}$ is the regulated output.

 The control objective is to design a low-frequency adaptive
controller to ensure that $y(t)$ tracks a given bounded continuous
reference signal $r(t)$ {\em both in transient and steady state},
while all other error signals remain bounded. More rigorously, the
control objective can be stated as  design of a low-pass control
signal $u(t)$ to achieve
\begin{equation}\label{DsDef}
y(s) \thickapprox D(s) r(s)\,,
\end{equation}
where $y(s), r(s)$ are Laplace transformations of $y(t), r(t)$
respectively, and $D(s)$ is a strictly proper stable LTI system
that specifies the desired transient and steady state performance.
We note that the control objective can be met if both the control
signal $u(t)$ and the system response $x(t)$ approximate the
corresponding signals of a LTI system with its response close to
$D(s)$.

\section{ $\mathcal{L}_1$ Adaptive Controller}\label{sec:Hinf}

 In this section, we
develop a novel adaptive control architecture that permits
complete transient characterization for system's both input and
output signals. Since $(A,b)$ is controllable, we choose $K$ to
ensure that $A_m  =  A - b K^{\top}$ is Hurwitz or, equivalently,
that
\begin{equation}\label{Hodef}
H_o(s) = (s I - A_m)^{-1} b
\end{equation}
is stable. The following control structure
\begin{eqnarray}\label{u1_design}
 u(t)  =  u_1(t)+u_2(t)\,, \quad
 u_1(t)  =  -K^{\top} x(t)\,,
\end{eqnarray}
where $u_2(t)$ is the adaptive controller to be determined later,
 leads to the following {\em
partially} closed-loop dynamics:
\begin{eqnarray}
\dot{x}(t) & = & A_m x(t)- b \theta^{\top} x(t) + b u_2(t), \quad x(0)=x_0 \label{problem_6}\\
y(t) & = & c^{\top} x(t)\,.\nonumber
\end{eqnarray}
For the linearly parameterized system in (\ref{problem_6}), we
consider the following state predictor
\begin{eqnarray}
\dot{\hat{x}}(t) & = & A_m \hat{x}(t)-b\hat{\theta}^{\top}(t)
x(t)+b u_2(t)\,,\quad \hat x(0)=x_0
\nonumber\\
\hat{y}(t) & = & c^{\top} \hat{x}(t)\label{L1_companionmodal}
\end{eqnarray}
along with the adaptive law for $\hat{\theta}(t)$:
\begin{equation}
\dot{\hat{\theta}}(t)  = \Gamma{\rm Proj}(\hat{\theta}(t),x(t)
\tilde{x}^{\top}(t) P b ), \quad {\hat{\theta}}(0)=\hat
\theta_0\,, \label{adaptivelaw_L1}
\end{equation}
where $\tilde{x}(t)=\hat x(t)-x(t)$ is the prediction error,
$\Gamma \in {\IR}^{n\times n}=\Gamma_c I_{n\times n}$ is the
matrix of adaptation gains, and $P$ is the solution of the
algebraic equation $ A_m^{\top} P+P A_m =- Q$, $Q>0$.

\begin{figure}[!h]
\begin{center}
\includegraphics[width=3.3in,height=2.0in]{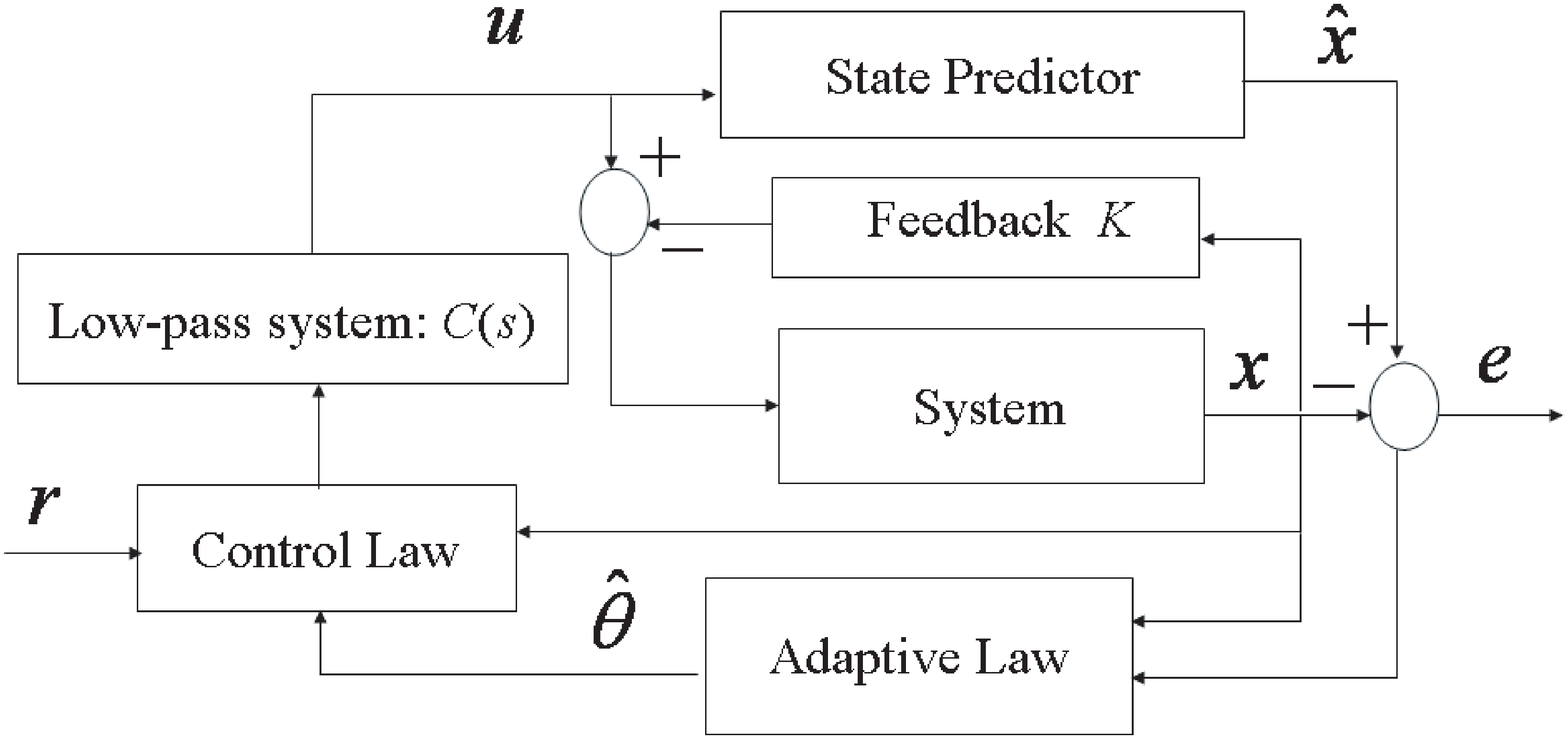}
\caption{Closed-loop system with $\mathcal{L}_1$ adaptive
controller}\label{fig:L1block}
\end{center}
\end{figure}

Letting
\begin{equation}\label{barrt_def}
\bar{r}(t)=\hat{\theta}^{\top}(t) x(t),
\end{equation}
the state predictor in (\ref{L1_companionmodal}) can be viewed as
a low-pass system with $u_2(t)$ being its control signal,
$\bar{r}(t)$ being a time-varying disturbance, which is not
prevented from having high-frequency oscillations.
We  consider the following control design for
(\ref{L1_companionmodal}):
\begin{equation}\label{us_def}
u_2(s)= C(s)\big( \bar{r}(s)+k_g r(s)\big) \,,
\end{equation}
where $u_2(s)$, $\bar{r}(s), r(s)$ are the Laplace transformations
of $u_2(t)$, $\bar{r}(t), r(t)$, respectively, $C(s)$ is a stable
and strictly proper system with low-pass gain $C(0)=1$, and $k_g$
is
\begin{equation}\label{kg_def}
 k_g = \lim_{s\rightarrow 0} \frac{1}{ c^{\top} H_o(s)} =\frac{1}{ c^{\top}
 H_o(0)}\,.
\end{equation}
The complete $\mathcal{L}_1$ adaptive controller consists of
(\ref{u1_design}), (\ref{L1_companionmodal}),
(\ref{adaptivelaw_L1}), (\ref{us_def}), and the closed-loop system
with it is illustrated in Fig. \ref{fig:L1block}.

Consider the closed-loop state predictor in
(\ref{L1_companionmodal}) with the control signal defined in
(\ref{us_def}). It can be viewed as an LTI system with two inputs
$r(t)$ and $\bar{r}(t)$:
\begin{eqnarray}
\hat{x}(s) & = & \bar{G}(s) \bar{r}(s)+G(s) r(s) \label{G1G2}\\
\bar{G}(s) & = & H_o(s) (C(s)-1) \label{G1s_def}\\
G(s) & = & k_g H_o(s) C(s)\,, \label{G2s_def}
\end{eqnarray}
where  $\hat{x}(s)$ is the Laplace transformation of $\hat{x}(t)$.
We note that $\bar{r}(t)$ is  related to $\hat{x}(t)$, $u(t)$ and
$r(t)$ via nonlinear relationships.
\begin{rem}\label{rem:stable}
Since both $H_o(s)$ and $C(s)$ are strictly proper stable systems,
one can check easily that $\bar{G}(s)$ and $G(s)$ are strictly
proper stable  systems, even though that $1-C(s)$ is proper.
\end{rem}
Let
\begin{equation}\label{thetabardef} \theta_{\max} =
\max_{\theta\in \Omega} \sum_{i=1}^{n} |\theta_i|\,,
\end{equation}
where $\theta_i$ is the $i^{th}$ element of $\theta$, $\Omega$ is
the compact set, where the unknown parameter lies. We now give the
$\mathcal{L}_1$ performance requirement that ensures
 stability of the entire system and desired transient performance, as discussed later in Section \ref{sec:convergence}.

{\bf $\mathcal{L}_1$-gain  requirement:} Design $K$ and $C(s)$ to
satisfy
\begin{equation}\label{condition3}
\Vert \bar{G}(s) \Vert_{\mathcal{L}_1}\theta_{\max}<1.
\end{equation}

\section{Analysis of $\mathcal{L}_1$ Adaptive
Controller}\label{sec:convergence}

\subsection{Stability and Asymptotic Convergence}

Consider the following Lyapunov function candidate:
\begin{equation}\label{lyapunov}
V(\tilde{x}(t),\tilde{\theta}(t))= \tilde{x}^{\top}(t) P
\tilde{x}(t)+ \tilde{\theta}^{\top}(t) \Gamma^{-1}
\tilde{\theta}(t)\,,
\end{equation}
where $P$ and $\Gamma$ are introduced in (\ref{adaptivelaw_L1}).
It follows from (\ref{problem_6}) and (\ref{L1_companionmodal})
that
\begin{equation}\label{L1errordynamics}
\dot{\tilde{x}}(t) = A_m \tilde{x}(t) - b \tilde{\theta}^{\top}(t)
x(t)\,, \quad \tilde x(0)=0\,.
\end{equation}
Hence, it is straightforward to verify from (\ref{adaptivelaw_L1})
that
\begin{equation}\label{Vdot}
\dot{V}(t) \leq -\tilde{x}^{\top}(t) Q \tilde{x}(t) \leq 0\,.
\end{equation}
 Notice that  the result
in (\ref{Vdot}) is independent of $u_2(t)$, however, one cannot
deduce stability from it. One needs to prove in addition that with
the $\mathcal{L}_1$ adaptive controller the state of the
predictor will remain bounded. Boundedness of the system state
then will follow.

\begin{thm}\label{thm1}
Given the system in (\ref{problem}) and the $\mathcal{L}_1$
adaptive controller defined via (\ref{u1_design}),
(\ref{L1_companionmodal}), (\ref{adaptivelaw_L1}), (\ref{us_def})
subject to (\ref{condition3}), the tracking error $\tilde x(t)$
converges to zero asymptotically:
\begin{equation}
\lim_{t\rightarrow \infty} \tilde{x}(t)=0.\label{thm_res}
\end{equation}
\end{thm}
{\bf Proof.} Let $\lambda_{\min}(P)$  be the minimum eigenvalue of
$P$. From (\ref{lyapunov}) and (\ref{Vdot}) it follows that $
\lambda_{\min}(P) \Vert \tilde{x}(t) \Vert^2\le
\tilde{x}^{\top}(t) P \tilde{x}(t) \leq V(t) \leq V(0), $ implying
that
\begin{equation}\label{thm_6}
\Vert \tilde{x}(t) \Vert^2 \leq V(0)/\lambda_{\min}(P),\;t\geq 0.
\end{equation}
From Definition \ref{defn1}, $ \displaystyle{\Vert \tilde{x}
\Vert_{{\mathcal L}_\infty} = \max_{i=1,..,n,t\geq 0}
|\tilde{x}_i(t)| }$. The relationship in (\ref{thm_6}) ensures
 $ \displaystyle{\max_{i=1,..,n,t\geq 0}} |\tilde{x}_i(t)|
\leq \sqrt{V(0)/\lambda_{\min}(P)}
$, 
and therefore for all $t>0$ one has $ \displaystyle{ \Vert
\tilde{x}_t \Vert_{{\mathcal L}_\infty} \leq
\sqrt{V(0)/\lambda_{\min}(P)}}\,.
$ 
Using the triangular relationship for norms implies that
\begin{equation}\label{thm_8}
|\;\Vert \hat{x}_t\Vert_{{\mathcal L}_\infty}-\Vert
x_t\Vert_{{\mathcal L}_\infty}\;|\leq
\sqrt{V(0)/\lambda_{\min}(P)}.
\end{equation}
The projection algorithm in (\ref{adaptivelaw_L1})  ensures that $
\hat{\theta}(t)\in \Omega,\forall t\geq 0 $. The definition of
$\bar r(t)$ in  (\ref{barrt_def}) implies that $ \Vert
\bar{r}_t\Vert_{{\mathcal L}_\infty} \leq \theta_{\max} \Vert
x_t\Vert_{{\mathcal L}_\infty}.
$ 
Substituting for  $\Vert x_t\Vert_{{\mathcal L}_\infty}$ from
(\ref{thm_8}) leads to the following
\begin{equation}\label{thm_10}
\Vert \bar{r}_t\Vert_{{\mathcal L}_\infty}\leq \theta_{\max}
\left(\Vert \hat{x}_t\Vert_{{\mathcal
L}_\infty}+\sqrt{V(0)/\lambda_{\min}(P)}\right).
\end{equation}
It follows from Lemma \ref{lem:L1} that
$ 
\Vert \hat{x}_t\Vert_{{\mathcal L}_\infty} \leq \Vert
\bar{G}(s)\Vert_{\mathcal{L}_1} \Vert \bar{r}_t\Vert_{{\mathcal
L}_\infty}+\Vert G(s)\Vert_{\mathcal{L}_1} \Vert
r_t\Vert_{{\mathcal L}_\infty}\,,
$ 
 which along with (\ref{thm_10}) gives the following upper bound
\begin{eqnarray}\label{thm_12}
&&\Vert \hat{x}_t\Vert_{{\mathcal L}_\infty} \leq \Vert
\bar{G}(s)\Vert_{\mathcal{L}_1} \theta_{\max} \left(\Vert
\hat{x}_t\Vert_{{\mathcal
L}_\infty}+\sqrt{V(0)/\lambda_{\min}(P)}\right)+\Vert
G(s)\Vert_{\mathcal{L}_1} \Vert r_t\Vert_{{\mathcal L}_\infty}.
\end{eqnarray}
Let
\begin{equation}\label{lambda}
\lambda=\Vert \bar{G}(s) \Vert_{\mathcal{L}_1}\theta_{\max}.
\end{equation}
 From (\ref{condition3}) it follows that
$\lambda<1$.
 The relationship in (\ref{thm_12}) can be written as
$ 
(1-\lambda)\Vert \hat{x}_t\Vert_{{\mathcal L}_\infty} \leq
\lambda\sqrt{V(0)/\lambda_{\min}(P)}+ \Vert
G(s)\Vert_{\mathcal{L}_1} \Vert r_t\Vert_{{\mathcal L}_\infty},
$ 
and hence
\begin{equation}\label{thm_15}
\Vert \hat{x}_t\Vert_{{\mathcal L}_\infty} \leq
(\lambda\sqrt{V(0)/\lambda_{\min}(P)}+\Vert
G(s)\Vert_{\mathcal{L}_1} \Vert r_t\Vert_{{\mathcal
L}_\infty})/(1-\lambda).
\end{equation}
Since $V(0), \lambda_{\min}(P), \Vert G(s) \Vert_{\mathcal{L}_1},
\Vert r\Vert_{{\mathcal L}_\infty}, \lambda$ are all finite and
$\lambda<1$, the relationship in (\ref{thm_15}) implies that
$\Vert \hat{x}_t\Vert_{{\mathcal L}_\infty}$ is finite for any
$t>0$, and hence $\hat{x}(t)$ is bounded. The relationship in
(\ref{thm_8}) implies that $ \Vert x_t\Vert_{{\mathcal L}_\infty}
$ is also finite for all $t>0$, and therefore $x(t)$ is bounded.
The adaptive law in (\ref{adaptivelaw_L1}) ensures that the
estimates $\hat{\theta}(t)$ are also bounded. From
(\ref{L1errordynamics}) it follows that  $\dot{\tilde{x}}(t)$ is
bounded, and it follows from Barbalat's lemma that $\displaystyle{
\lim_{t\rightarrow \infty} \tilde{x}(t) =0} $. $\hfill{\square}$

%
%

\subsection{Reference System}

\begin{figure}[!t]
\begin{center}
\includegraphics[width=3.4in,height=2.0in]{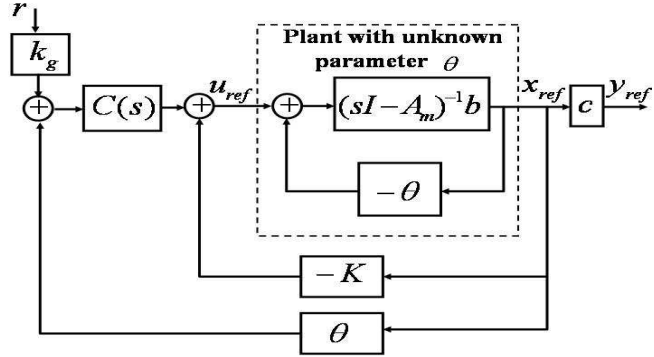}
\caption{Closed-loop  reference LTI system }\label{fig:refblock}
\end{center}
\end{figure}

In this section we characterize the reference system that the
$\mathcal{L}_1$ adaptive controller in (\ref{u1_design}),
(\ref{L1_companionmodal}), (\ref{adaptivelaw_L1}), (\ref{us_def})
tracks both in transient and steady state, and this tracking is
valid for system's both input and output signals. Towards that
end, consider the following {\em ideal} version of the adaptive
controller in (\ref{u1_design}), (\ref{us_def}):
\begin{equation}\label{ref_u}
u_{ref}(s) = C(s) \left( k_g r(s) + \eta(s)\right) - K^{\top}
x_{ref}(s)\,,
\end{equation}
where $\eta(s)$ is the Laplace transformation of
\begin{equation}\label{eta_def}
\eta(t)=\theta^{\top} x_{ref}(t),
\end{equation}
 and $x_{ref}(s)$ is used to
denote  the Laplace transformation of the state  $x_{ref}(t)$ of
the closed-loop system. The closed-loop system (\ref{problem})
 with the controller (\ref{ref_u}) is given in  Fig. \ref{fig:refblock}.

\begin{rem}
Notice that when $C(s)=1$, one recovers the reference model of
MRAC.
If $C(s)\neq 1$, then the control law in (\ref{ref_u}) changes the
bandwidth of the ideal control signal
$u_{ideal}(t)=-K^\top x(t)+\theta^{\top} x(t)+k_g r(t)$. 
\end{rem}
The control law in (\ref{ref_u}) leads to the following
closed-loop dynamics:
\begin{eqnarray}
x_{ref}(s) & = & H_o(s) \left(k_g C(s) r(s)+ (C(s)-1)
\theta^{\top} x_{ref}(s)\right)\label{yrefdef}\\
y_{ref}(s) & = &  c^{\top}  x_{ref}(s)\,,\nonumber
\end{eqnarray}
which  can be explicitly solved for $ x_{ref}(s) = \Big( I -
(C(s)-1) H_o(s) \theta^{\top} \Big)^{-1} H_o(s) k_g C(s) r(s). $
Hence, it follows from (\ref{G1s_def}) and (\ref{G2s_def})
 that
\begin{equation}\label{ref_x}
x_{ref}(s) = (I-\bar{G}(s)\theta^{\top})^{-1} G(s) r(s)\,.
\end{equation}

\begin{lem}\label{lem:refstable}
If $\Vert \bar{G}(s) \Vert_{\mathcal{L}_1}\theta_{\max}<1$, then
\begin{eqnarray}
(i) & (I-\bar{G}(s)\theta^{\top})^{-1} \;\;\hbox{ is stable};
\label{lemrefstable_0}\\
(ii) & (I-\bar{G}(s)\theta^{\top})^{-1}G(s) \;\;\hbox{ is
stable}.\nonumber 
\end{eqnarray}
\end{lem}
{\bf Proof.}
%
It follows from (\ref{L1def}) that
$$ 
\Vert \bar{G}(s)\theta^{\top} \Vert_{\mathcal{L}_1} =
\max_{i=1,..,n} \left(
 \Vert \bar{G}_i(s)\Vert_{\mathcal{L}_1} \left(\sum_{j=1}^{n}
 |\theta_j|\right)\right)\,,
$$ 
where $ \bar{G}_i(s) $ is the $i^{th}$ element of $G(s)$, and
$\theta_j$ is the $j^{th}$ element of $\theta$. From
(\ref{thetabardef}) we have $ \sum_{j=1}^{n} |\theta_j| \leq
\theta_{\max}\,, $ and hence
\begin{equation}\label{lemrefstable_12}
\Vert \bar{G}(s)\theta^{\top} \Vert_{\mathcal{L}_1} \leq
\max_{i=1,..,n} \Big(
 \Vert \bar{G}_i(s)\Vert_{\mathcal{L}_1} \Big) \theta_{\max} = \Vert \bar{G}(s)
 \Vert_{\mathcal{L}_1} \theta_{\max}, \quad \forall ~\theta\in \Omega.
\end{equation}
The relationship in (\ref{condition3}) implies that
$ 
\Vert \bar{G}(s) \theta^{\top} \Vert_{\mathcal{L}_1} <  1,
$ 
and therefore Theorem \ref{thm:sg} ensures
 that the LTI system
$(I-\bar{G}(s)\theta^{\top})^{-1}$ is stable. Since $G(s)$ is
stable, Remark \ref{rem:stable} implies that
$(I-\bar{G}(s)\theta^{\top})^{-1}G(s)$ is stable.
$\hfill{\square}$

\subsection{System Response and Control Signal of the $\mathcal{L}_1$ Adaptive Controller}

Letting
\begin{equation}\label{r1_def}
r_1(t)= \tilde{\theta}^{\top}(t) x(t),
\end{equation}
we notice that $\bar r(t)$ in (\ref{barrt_def}) can be rewritten
as $ \bar{r}(t)=\theta^{\top}(\hat{x}(t)-\tilde{x}(t))+r_1(t)$.
Hence, the state predictor in (\ref{G1G2}) can be rewritten as $
\hat{x}(s) = \bar{G}(s)\left(\theta^{\top}
\hat{x}(s)-\theta^{\top} \tilde{x}(s)+r_1(s)\right)+G(s) r(s) $,
where $r_1(s)$ is the Laplace transformation of $r_1(t)$ defined
in (\ref{r1_def}), and further put into the form:
\begin{equation}\label{lem2_3}
\hat{x}(s) = (I-\bar{G}(s)\theta^{\top})^{-1}\Big(
-\bar{G}(s)\theta^{\top} \tilde{x}(s)+\bar{G}(s)r_1(s)+G(s)
r(s)\Big).
\end{equation}
It follows from (\ref{L1errordynamics}) and (\ref{r1_def}) that $
\dot{\tilde{x}}(t) = A_m \tilde{x}(t) - b r_1(t)\,, $ and hence
\begin{equation}\label{lem2_5}
\tilde{x}(s) = - H_o(s) r_1(s).
\end{equation}
Using the expression of $\bar G(s)$ from (\ref{G1s_def}), the
state of the  predictor can be presented as
$
\hat{x}(s) = (I-\bar{G}(s)\theta^{\top})^{-1}G(s) r(s)
+(I-\bar{G}(s)\theta^{\top})^{-1} \left( -\bar{G}(s)\theta^{\top}
\tilde{x}(s)-(C(s)-1) \tilde{x}(s)\right) .
$
Using $x_{ref}(s)$ from (\ref{ref_x}) and recalling the definition
of $\tilde{x}(s) = \hat{x}(s)-x(s)$, one arrives at
\begin{equation}\label{lem2_6}
x(s) = x_{ref}(s) -\Big(I+(I-\bar{G}(s)\theta^{\top})^{-1} \big(
\bar{G}(s)\theta^{\top}+ (C(s)-1) I  \big)\Big) \tilde{x}(s).
\end{equation}
The expressions in (\ref{u1_design}), (\ref{us_def}) and
(\ref{ref_u}) lead to the following expression of the control
signal
\begin{equation}\label{lem2_8}
u(s) = u_{ref}(s) + C(s) r_1(s)+(C(s) \theta^{\top}-K^{\top})
(x(s)-x_{ref}(s))\,.
\end{equation}

\subsection{Asymptotic Performance and Steady State Error}

\begin{thm}\label{thm:3}
Given the system in (\ref{problem}) and the $\mathcal{L}_1$
adaptive controller defined via  (\ref{u1_design}),
(\ref{L1_companionmodal}), (\ref{adaptivelaw_L1}), (\ref{us_def})
 subject to (\ref{condition3}), we have:
\begin{eqnarray}
\lim_{t\rightarrow \infty} \Vert x(t)-x_{ref}(t) \Vert & = & 0\,,
\label{thm3_01} \\
\lim_{t\rightarrow \infty} \vert u(t)-u_{ref}(t) \vert & = & 0\,.
\label{thm3_02}
\end{eqnarray}
\end{thm}
{\bf Proof.} Let
\begin{equation}\label{thm3_3new}
r_2(t)=x_{ref}(t)-x(t)\,.
\end{equation}
It follows from (\ref{lem2_6}) that
\begin{equation}\label{r2sdef}
r_2(s) =\left(I+(I-\bar{G}(s)\theta^{\top})^{-1} \left(
\bar{G}(s)\theta^{\top}+ (C(s)-1) I  \right)\right) \tilde{x}(s).
\end{equation}
The signal $r_2(t)$ can be viewed as the  response of the LTI
system
\begin{equation}\label{H2sdef}
H_2(s)=I+(I-\bar{G}(s)\theta^{\top})^{-1} \left(
\bar{G}(s)\theta^{\top}+ (C(s)-1) I  \right)
\end{equation}
 to the bounded error
signal $\tilde{x}(t)$. It follows from (\ref{lemrefstable_0}) and
Remark \ref{rem:stable} that $(I-\bar{G}(s)\theta^{\top})^{-1}$,
$\bar{G}(s)$, $C(s)$ are stable and, therefore, $H_2(s)$ is
stable. Hence,  from (\ref{thm_res}) we have $
\displaystyle{\lim_{t \rightarrow \infty} r_2(t)=0}$.
Let
\begin{equation}\label{thm3_5}
r_3(s) =C(s) r_1(s)+(C(s) \theta^{\top}-K^{\top})
(x(s)-x_{ref}(s)).
\end{equation}
It follows from (\ref{lem2_8}) that
\begin{equation}\label{thm3_6}
r_3(t)=u(t)-u_{ref}(t).
\end{equation}
Since the projection operator ensures that $\tilde{\theta}(t)$ is
bounded, it follows from (\ref{L1errordynamics}) and
(\ref{thm_res})  that $ \displaystyle{\lim_{t \rightarrow \infty}
r_1(t) =0.} $ Since $C(s)$ is a stable proper system, it follows
from (\ref{thm3_01}) that $ \displaystyle{\lim_{t \rightarrow
\infty}} r_3(t) =0$. $\hfill{\square}$

\begin{lem}\label{lem:5}
Given the system in (\ref{problem}) and the $\mathcal{L}_1$
adaptive controller defined via (\ref{u1_design}),
(\ref{L1_companionmodal}), (\ref{adaptivelaw_L1}), (\ref{us_def})
 subject to (\ref{condition3}), if
$r(t)$ is constant, then $ \displaystyle{ \lim_{t\rightarrow
\infty} y(t)=r. }$
\end{lem}
{\bf Proof.} Since
$
y_{ref}(t)= c^{\top} x_{ref}(t),
$
 it follows from (\ref{thm3_01})  that
$ \displaystyle{\lim_{t \rightarrow \infty} (y(t)- y_{ref}(t))=0.}
$ From (\ref{ref_x}) it follows that
$ 
y_{ref}(s) = c^{\top} (I-\bar{G}(s)\theta^{\top})^{-1}G(s) r(s).
$ 
The end value theorem ensures
\begin{equation}\label{lem5_8}
\lim_{t\rightarrow \infty} y_{ref}(t) = \lim_{s\rightarrow 0}
c^{\top} (I-\bar{G}(s)\theta^{\top})^{-1}G(s)  r = c^{\top} H_o(0)
C(0) k_g r.
\end{equation}
Definition of $k_g$ in (\ref{kg_def}) leads to $
\displaystyle{\lim_{t \rightarrow \infty} y(t)=r. }$
$\hfill{\square}$

\subsection{Bounded Tracking Error Signal}

\begin{lem}\label{lem:1}
Let  $\Gamma=\Gamma_c \II$, where $\Gamma_c\in {\IR}^{+}$, and
$\II$ is the identity matrix. For the system in (\ref{problem})
\begin{equation}
|| \tilde{x}(t) ||  \leq
\sqrt{\bar{\theta}_{\max}/(\lambda_{\min}(P) \Gamma_c)}\, ,\quad
\bar{\theta}_{\max} \triangleq \max_{\theta\in \Omega}
\sum_{i=1}^{n} 4 \theta_i^2,\qquad \forall t\geq 0,
\label{barthemax}
\end{equation}
and $\lambda_{\min}(P)$ is the minimum eigenvalue of $P$.
\end{lem}
{\bf Proof.}
For the $V(t)$  in (\ref{lyapunov}), the following upper bound is
straightforward to derive: $ \tilde{x}^{\top}(t) P \tilde{x}(t)
\leq  V(t) \leq V(0), \forall t \geq 0. $ The projection algorithm
ensures that $ \hat{\theta}(t) \in \Omega , ~ \forall t \geq 0, $
and therefore
\begin{equation}\label{lem1_5}
\max_{t\geq 0} \tilde{\theta}^{\top}(t) \Gamma^{-1}
\tilde{\theta}(t) \leq \frac{\bar{\theta}_{\max}}{\Gamma_c },
\quad \forall t\geq 0,
\end{equation}
where $\bar{\theta}_{\max}$ is defined in (\ref{barthemax}). Since
$\tilde x(0)=0$, then $ V(0) = \tilde{\theta}^{\top}(0)
\Gamma^{-1} \tilde{\theta}(0)$, which leads to
$ 
\displaystyle{\tilde{x}^{\top}(t) P \tilde{x}(t) \leq
\frac{\bar{\theta}_{\max}}{\Gamma_c },~ t\geq 0.}
$ 
Since $ \lambda_{\min}(P) \Vert \tilde{x} \Vert^2 \leq
\tilde{x}^{\top}(t) P \tilde{x}(t), $ then
$ 
\displaystyle{|| \tilde{x}(t) ||  \leq
\sqrt{\frac{\bar{\theta}_{\max}}{\lambda_{\min}(P) \Gamma_c}}.}
$ 

\subsection{Transient Performance}\label{sec:ana_tp}

We note that $(A-b K^{\top}, b)$ is the state space realization of
$H_o(s)$. Since $(A, b)$ is controllable, it can be easily proved
 that $(A-b K^{\top}, b)$ is also controllable. It follows
from Lemma \ref{lem:pre2} that there exists $c_{o}\in {\IR}^n$
such that
\begin{equation}\label{thm5_pre}
c_{o}^{\top} H_o(s) = N_n(s)/N_d(s)\,,
\end{equation}
where the order of $N_d(s)$ is one more than the order of
$N_n(s)$, and both $N_n(s)$ and $N_d(s)$ are stable polynomials.

\begin{thm}\label{thm:5}
Given the system in (\ref{problem}) and the $\mathcal{L}_1$
adaptive controller defined via (\ref{u1_design}),
(\ref{L1_companionmodal}), (\ref{adaptivelaw_L1}), (\ref{us_def})
subject to (\ref{condition3}), we have:
\begin{eqnarray}
\Vert x-x_{ref} \Vert_{{\mathcal L}_{\infty}}  & \leq &
\gamma_1/\sqrt{\Gamma_c}\, ,\label{thm5_00} \\
\Vert y-y_{ref} \Vert_{{\mathcal L}_{\infty}}  & \leq &
\Vert c^{\top}\Vert_{{\mathcal L}_1} \gamma_1/\sqrt{\Gamma_c}\, ,\label{thm5_02}\\
 \Vert u -
u_{ref} \Vert_{{\mathcal L}_{\infty}} & \leq &
\gamma_2/\sqrt{\Gamma_c}\, ,\label{thm5_01}
\end{eqnarray}
where $\Vert c^{\top}\Vert_{{\mathcal L}_1}$ is the
$\mathcal{L}_1$ gain of $c^{\top}$ and
\begin{eqnarray}
\gamma_1 & = &
\|H_2(s)\|_{\mathcal L_1}\sqrt{\bar{\theta}_{\max}/\lambda_{\max}(P)}\,, \label{gamma1def}\\
\gamma_2 & = & \Big\|C(s) \frac{1}{c_{o}^{\top} H_o(s)}
c_{o}^{\top}
\Big\|_{\mathcal{L}_1}\sqrt{\bar{\theta}_{\max}/\lambda_{\max}(P)}+\|C(s)
\theta^{\top} -K^{\top} \|_{\mathcal{L}_1} \gamma_1\,.
\label{gamma2def}
\end{eqnarray}
\end{thm}
\medskip

{\bf Proof.} It follows from (\ref{r2sdef}), (\ref{H2sdef}) and
Lemma \ref{lem:L1_ext} that
$ 
\Vert r_2\Vert_{{\mathcal L}_{\infty}} \leq \Vert H_2(s)
\Vert_{\mathcal{L}_1} \Vert \tilde{x} \Vert_{{\mathcal
L}_{\infty}},
$ 
while Lemma \ref{lem:1} implies that
\begin{equation}\label{boundxtilde}
\Vert \tilde{x} \Vert_{{\mathcal L}_{\infty}} \leq
\sqrt{\bar{\theta}_{\max}/(\lambda_{\max}(P) \Gamma_c)}\,.
\end{equation}
Therefore,
$ 
\Vert r_2\Vert_{{\mathcal L}_{\infty}} \leq  \Vert H_2(s)
\Vert_{\mathcal{L}_1}
\sqrt{\frac{\bar{\theta}_{\max}}{\lambda_{\max}(P) \Gamma_c}},
$ 
which   leads to (\ref{thm5_00}). The upper bound in
(\ref{thm5_02}) follows from (\ref{thm5_00}) and Lemma
\ref{lem:L1cas} directly. From (\ref{lem2_5}) we have
\begin{eqnarray}
& & r_3(s)  = C(s) \frac{1}{c_{o}^{\top} H_o(s)} c_{o}^{\top}
H_o(s)
r_1(s)     + (C(s) \theta^{\top}-K^{\top})(x(s)-x_{ref}(s)) \nonumber\\
& & \qquad =  -C(s) \frac{1}{c_{o}^{\top} H_o(s)} c_{o}^{\top}
\tilde{x}(s)+ (C(s) \theta^{\top}
-K^{\top})(x(s)-x_{ref}(s))\,,\nonumber
\end{eqnarray}
where $c_{o}$ is introduced in (\ref{thm5_pre}).
 It follows from (\ref{thm5_pre}) that
$ 
 C(s) \frac{1}{c_{o}^{\top} H_o(s)} = C(s)
 \frac{N_d(s)}{N_n(s)}\,,
$ 
where $N_d(s)$, $N_n(s)$ are stable polynomials and the order of
$N_n(s)$ is one less than the order of $N_d(s)$.  Since $C(s)$ is
stable and strictly proper, the complete system
$
 C(s) \frac{1}{c_{o}^{\top} H_o(s)}
$ is proper and stable, which implies that its $\mathcal{L}_1$
gain exists and is finite. Hence, we have
$
\Vert r_3\Vert_{{\mathcal L}_{\infty}} \leq \Big\| C(s)
\frac{1}{c_{o}^{\top} H_o(s)} c_{o}^{\top} \Big\|_{\mathcal{L}_1}
\Vert \tilde{x} \Vert_{{\mathcal L}_{\infty}}+ \Vert C(s)
\theta^{\top} -K^{\top}\Vert_{{\mathcal L}_1} \Vert x-x_{ref}
\Vert_{{\mathcal L}_{\infty}}\,.
$
 Lemma \ref{lem:1}  leads to the upper bound in (\ref{thm5_01}):
$$ 
\Vert r_3\Vert_{{\mathcal L}_{\infty}} \leq  \Big\| C(s)
\frac{1}{c_{o}^{\top} H_o(s)} c_{o}^{\top} \Big\|_{\mathcal{L}_1}
\sqrt{\frac{\bar{\theta}_{\max}}{\lambda_{\max}(P) \Gamma_c}}+
\Vert C(s) \theta^{\top} -K^{\top} \Vert_{{\mathcal L}_1} \Vert
x-x_{ref} \Vert_{{\mathcal L}_{\infty}}\,.
$$ 
 $\hfill{\square}$

\begin{cor}\label{cor:1}
Given the system in (\ref{problem}) and the $\mathcal{L}_1$
adaptive controller defined via (\ref{u1_design}),
(\ref{L1_companionmodal}), (\ref{adaptivelaw_L1}), (\ref{us_def})
subject to (\ref{condition3}), we have:
\begin{eqnarray}
\lim_{\Gamma_c\rightarrow \infty}  \left(x(t)-x_{ref}(t) \right)&
= & 0\, ,\qquad \forall t\geq 0,\label{cor1_00}
\\
\lim_{\Gamma_c\rightarrow \infty}  \left(y(t)-y_{ref}(t) \right)&
= & 0\, ,\qquad \forall t\geq 0,\label{cor1_02}
\\
\lim_{\Gamma_c\rightarrow \infty} \left(u(t)-u_{ref}(t)\right) & =
& 0\,,\qquad \forall t\geq 0\,.\label{cor1_01}
\end{eqnarray}
\end{cor}

%
Corollary \ref{cor:1} states that $x(t)$, $y(t)$ and $u(t)$ follow
$x_{ref}(t)$, $y_{ref}(t)$ and $u_{ref}(t)$ not only
asymptotically but also during the transient, provided that the
adaptive gain is selected sufficiently large. Thus, the control
objective is reduced to designing $K$ and $C(s)$ to ensure that
the reference LTI system has the desired  response $D(s)$.

\begin{rem}
Notice that if we set $C(s)=1$, then the $\mathcal{L}_1$ adaptive
controller degenerates into a MRAC type. In that case $\Big\|C(s)
\frac{1}{c_{o}^{\top} H_o(s)} c_{o}^{\top} \Big\|_{\mathcal{L}_1}
$ cannot be finite, since $H_o(s)$ is strictly proper. Therefore,
from (\ref{gamma2def}) it follows that $\gamma_2 \rightarrow
\infty$, and hence for the control signal in  MRAC one can not
reduce the bound in (\ref{thm5_01}) by increasing the adaptive
gain.
\end{rem}

\section{Design of the $\mathcal{L}_1$ Adaptive Controller}\label{sec:design}

We proved that the error between the state and the control signal
of the closed-loop system with $\mathcal{L}_1$ adaptive controller
in (\ref{problem}), (\ref{u1_design}), (\ref{L1_companionmodal}),
(\ref{adaptivelaw_L1}), (\ref{us_def}) (Fig. \ref{fig:L1block})
and  the state and the control signal of the closed-loop reference
system in (\ref{ref_u}), (\ref{ref_x}) (Fig. \ref{fig:refblock})
can be rendered arbitrarily small by choosing large adaptive gain.
Therefore, the control objective is reduced to determining $K$ and
$C(s)$ to ensure that the reference system in (\ref{ref_u}),
(\ref{ref_x}) (Fig. \ref{fig:refblock}) has the desired response
$D(s)$
 from $r(t)$ to $y_{ref}(t)$. Notice that the
reference system in Fig. \ref{fig:refblock} depends upon the
unknown parameter $\theta$.

Consider the following signals:
\begin{eqnarray}
y_{des}(s) &=& c^{\top} G(s) r(s) = C(s) k_g c^{\top} H_o(s)
r(s)\,,\label{ydes_def}\\
u_{des}(s) &=&k_g C(s) \left( 1  + C(s) \theta^{\top} H_o(s)
-K^{\top} H_o(s)\right)r(s)\label{udes_def}.
\end{eqnarray}
We note that $u_{des}(t)$ depends on the unknown parameter
$\theta$, while $y_{des}(t)$ does not.

\begin{lem}\label{lem:design}
Subject to (\ref{condition3}), the following upper bounds hold:
\begin{eqnarray}
\Vert y_{ref} - y_{des} \Vert_{\mathcal{L}_\infty}  & \leq &
\frac{\lambda}{1-\lambda} \Vert c^{\top} \Vert_{\mathcal{L}_1}
\Vert G(s) \Vert_{\mathcal{L}_1}\Vert r \Vert_{\mathcal{L}_\infty}
,\label{lemdesign_00}
\\
\Vert y_{ref} - y_{des} \Vert_{\mathcal{L}_\infty}  & \leq &
\frac{1}{1-\lambda} \Vert c^{\top} \Vert_{\mathcal{L}_1} \Vert h_3
\Vert_{\mathcal{L}_\infty}\, , \label{lemdesign_01} \\
\Vert u_{ref} - u_{des} \Vert_{\mathcal{L}_\infty}  & \leq &
\frac{\lambda}{1-\lambda} \Vert C(s) \theta^{\top}
-K^{\top}\Vert_{\mathcal{L}_1} \Vert G(s)
\Vert_{\mathcal{L}_1}\Vert r \Vert_{\mathcal{L}_\infty}\,,
\label{lemdesign_03} \\
\Vert u_{ref} - u_{des} \Vert_{\mathcal{L}_\infty}  & \leq &
\frac{1}{1-\lambda} \Vert C(s) \theta^{\top}
-K^{\top}\Vert_{\mathcal{L}_1}
 \Vert h_3 \Vert_{\mathcal{L}_\infty}\,,
\label{lemdesign_02}
\end{eqnarray}
 where $\lambda$ is defined in (\ref{lambda}),
and $h_3(t)$ is the inverse Laplace transformation of
\begin{equation}\label{H3sdef}
H_3(s) = (C(s)-1) C(s) r(s) k_g  H_o(s) \theta^{\top} H_o(s).
\end{equation}
\end{lem}
{\bf Proof.} It follows from (\ref{yrefdef}) and (\ref{ref_x})
that
$ 
y_{ref}(s) = c^{\top} (I-\bar{G}(s)\theta^{\top})^{-1} G(s) r(s).
$ 
Following Lemma \ref{lem:refstable}, the condition in
(\ref{condition3}) ensures the stability of the reference LTI
system. Since $(I-\bar{G}(s)\theta^{\top})^{-1}$ is stable, then
one can expand it into convergent series and further write
\begin{eqnarray}\label{lemdesign_3}
y_{ref}(s)  =  c^{\top} \left(I + \sum_{i=1}^{\infty}
(\bar{G}(s)\theta^{\top})^i \right) G(s) r(s) =  y_{des}(s) +
c^{\top} \left( \sum_{i=1}^{\infty} (\bar{G}(s)\theta^{\top})^i
\right) G(s) r(s).
\end{eqnarray}
Let $r_4(s)= c^{\top} \left( \sum_{i=1}^{\infty}
(\bar{G}(s)\theta^{\top})^i \right) G(s) r(s). $ Then
\begin{equation}\label{lemdesign_5}
r_4(t) = y_{ref}(t)- y_{des}(t) , ~~~~~\forall t\geq 0.
\end{equation}
The relationship in (\ref{lemrefstable_12})  implies that
$
\Vert \bar{G}(s)\theta^{\top} \Vert_{\mathcal{L}_1} \leq \lambda
$,
   and it follows from Lemma \ref{lem:L1cas} that
\begin{equation}\label{lemdesign_6}
\Vert r_4 \Vert_{\mathcal{L}_\infty} \leq
\left(\sum_{i=1}^{\infty} \lambda^i \right) \Vert c^{\top}
\Vert_{\mathcal{L}_1} \Vert G \Vert_{\mathcal{L}_1}\Vert r
\Vert_{\mathcal{L}_\infty}=\frac{\lambda}{1-\lambda} \Vert
c^{\top} \Vert_{\mathcal{L}_1} \Vert G \Vert_{\mathcal{L}_1}\Vert
r \Vert_{\mathcal{L}_\infty}.
\end{equation}
From (\ref{lemdesign_3})  we have $ y_{ref}(s)=y_{des}(s) +
c^{\top} \Big( \sum_{i=1}^{\infty} (\bar{G}(s)\theta^{\top})^{i-1}
\Big) \bar{G}(s)\theta^{\top} G(s) r(s), $ which along with
(\ref{G1s_def}), (\ref{G2s_def}) and (\ref{H3sdef}) leads to
\begin{eqnarray}
y_{ref}(s) =
  y_{des}(s) + c^{\top} \left(
\sum_{i=1}^{\infty} (\bar{G}(s)\theta^{\top})^{i-1} \right)
H_3(s). \nonumber
\end{eqnarray}
 Lemma \ref{lem:L1_ext} immediately implies that
$ \Vert r_4 \Vert_{\mathcal{L}_\infty} \leq
\Big(\sum_{i=1}^{\infty} \lambda^{i-1}\Big) \Vert c^{\top}
\Vert_{\mathcal{L}_1} \Vert h_3 \Vert_{\mathcal{L}_\infty}. $
Comparing $u_{des}(s)$ in (\ref{udes_def}) to $u_{ref}(s)$ in
(\ref{ref_u}) it follows that $u_{des}(s)$ can be written as $
u_{des}(s) = k_g C(s) r(s) + (C(s) \theta^{\top} -K^{\top})
x_{des}(s)\,, $ where $ x_{des}(s) = C(s) k_g  H_o(s) r(s). $
Therefore $ u_{ref}(s)-u_{des}(s) = (C(s) \theta^{\top} -K^{\top})
(x_{ref}(s)-x_{des}(s)). $ Hence, it follows from Lemma
\ref{lem:L1} that $ \Vert u_{ref}-u_{des}
\Vert_{\mathcal{L}_{\infty}} \leq \Vert C(s) \theta^{\top}
-K^{\top}\Vert_{\mathcal{L}_1} \Vert
x_{ref}-x_{des}\Vert_{\mathcal{L}_{\infty}}. $ Using the same
steps as for $\Vert y_{ref}-y_{des}\Vert_{\mathcal{L}_{\infty}}$,
we have
\begin{eqnarray}
\Vert x_{ref} - x_{des} \Vert_{\mathcal{L}_\infty}  & \leq &
\frac{\lambda}{1-\lambda} \Vert G(s) \Vert_{\mathcal{L}_1}\Vert r
\Vert_{\mathcal{L}_\infty} ,\nonumber
\\
\Vert x_{ref} - x_{des} \Vert_{\mathcal{L}_\infty}  & \leq &
\frac{1}{1-\lambda}  \Vert h_3 \Vert_{\mathcal{L}_\infty}\, ,
\nonumber
\end{eqnarray}
which leads to the upper bounds in (\ref{lemdesign_03}) and
(\ref{lemdesign_02}). $\hfill{\square}$

Thus, the problem is reduced to finding a strictly proper stable
$C(s)$ to ensure that
\begin{eqnarray}
(i)  & &   \lambda<1 \hbox{ or } \Vert h_3 \Vert_{\mathcal{L}_\infty} \hbox{ are sufficiently small}, \label{desre_5}\\
(ii) & &  y_{des}(s) \thickapprox D(s) r(s)\,, \label{desre_6}
\end{eqnarray}
where $D(s)$ is the desired LTI system introduced in
(\ref{DsDef}). Then, Theorem \ref{thm:5} and Lemma
\ref{lem:design} will imply that the output $y(t)$ of the system
in (\ref{problem}) and the $\mathcal{L}_1$ adaptive control signal
$u(t)$ will follow $y_{des}(t)$ and $u_{des}(t)$ both in transient
and steady state with quantifiable bounds, given in
(\ref{thm5_02}), (\ref{thm5_01}) and
(\ref{lemdesign_00})-(\ref{lemdesign_02}).

Notice that  $\lambda <1$ is required for stability. From
(\ref{ydes_def})-(\ref{lemdesign_02}), it follows that for
achieving $y_{des}(s) \thickapprox D(s) r(s)$ it is desirable to
ensure that $\lambda$ or $\Vert h_3 \Vert_{\mathcal{L}_\infty}$
are sufficiently small and, in addition, $C(s) c^{\top} H_o(s)
\thickapprox D(s)$. We notice that these requirements are not in
conflict with each other.  So, using Lemma
 \ref{lem:L1cas}, one can consider the following conservative upper bound
\begin{equation}\label{design_md1}
\lambda = \Vert \bar{G}(s)\Vert_{\mathcal{L}_1} \theta_{\max}=
\Vert
 H_o(s)
( C(s)-1)\Vert_{\mathcal{L}_1} \theta_{\max}\leq \Vert H_o(s)
\Vert_{\mathcal{L}_1} \Vert C(s)-1 \Vert_{\mathcal{L}_1}
\theta_{\max}.
\end{equation}
Thus,  minimization of $\lambda$ can be  achieved from two
different perspectives: i) fix $C(s)$ and minimize $\Vert H_o(s)
\Vert_{\mathcal{L}_1}$, ii) fix $H_o(s)$ and minimize the
$\mathcal L_1$-gain of one of the cascaded systems $ \|H_o(s)(
C(s)-1) \|_{\mathcal{L}_1}$, $\Vert (C(s)-1) r(s)
\Vert_{\mathcal{L}_1}$ or $ \Vert C(s)
(C(s)-1)\Vert_{\mathcal{L}_1}$ via the choice of $C(s)$.

{\em i) High-gain design.} Set $C(s)=D(s)$. Then minimization of
$\Vert H_o(s) \Vert_{\mathcal{L}_1}$ can be achieved via high-gain
feedback by choosing $K$ sufficiently large.   However, minimized
$\Vert H_o(s) \Vert_{\mathcal{L}_1} $ via large $K$ leads to
high-gain design with reduced robustness properties. Since $C(s)$
is a strictly proper system containing the dominant poles of the
closed-loop system in $ k_g c^{\top} H_o(s) C(s)$ and $k_g
c^{\top} H_o(0)=1$, we have
$ 
 k_g c^{\top} H_o(s) C(s) \approx C(s) =D(s).
$ 
Hence, the system response will be $ y_{ref}(s) \approx D(s) r(s).
$ We note that with large feedback $K$, the performance of
$\mathcal{L}_1$ adaptive controller degenerates into a high-gain
type. The shortcoming of this design is that the high gain
feedback $K$ leads to a reduced phase and time-delay margin and
consequently affects robustness.

{\em ii) Design without high-gain feedback. } As in  MRAC,  assume
that we can select $K$ to ensure
\begin{equation}\label{L1_app2}
k_g c^{\top} H_o(s) \thickapprox D(s).
\end{equation}


\begin{lem}\label{lem:designgain}
Let
\begin{equation}\label{firstorderCs}
C(s) =\frac{\omega}{s+\omega}\,.
\end{equation}
For any single input $n$-output strictly proper stable system
$H_o(s)$ the following is true: $$ \displaystyle{
\lim_{\omega\rightarrow \infty} \Vert ( C(s) -1) H_o(s)
\Vert_{\mathcal{L}_1} =0.} $$
\end{lem}
{\bf Proof.} It follows from (\ref{firstorderCs}) that
$\displaystyle{ (C(s)-1) H_o(s) = \frac{-s}{s+\omega} H_o(s) =
\frac{-1}{s+\omega} s H_o(s).}
$ 
Since $H_o(s)$ is strictly proper and stable, $s H_o(s)$ is stable
and has relative degree $\geq 0$, and hence $\Vert s H_o(s)
\Vert_{\mathcal{L}_1}$ is finite. Since $ \displaystyle{ \Big\|
\frac{-1}{s+\omega} \Big\|_{\mathcal{L}_1} = \frac{1}{\omega}\,,}
$ it follows from (\ref{lem:L1cas}) that $ \displaystyle{ \Vert
(C(s)-1) H_o(s) \Vert_{\mathcal{L}_1} \leq \frac{1}{\omega} \Vert
s H_o(s) \Vert_{\mathcal{L}_1}} $, and  the proof is complete.
$\hfill{\square}$

Lemma \ref{lem:designgain} states that if one chooses $k_g
c^{\top} H_o(s) r(s) \thickapprox D(s)$, then by increasing the
bandwidth of the low-pass system $C(s)$, it is possible to render
$ \Vert \bar{G}(s) \Vert_{\mathcal{L}_1}$ arbitrarily small.
 With large $\omega$, the pole $-\omega$ due to $C(s)$ is
omitted, and $H_o(s)$ is the dominant reference system leading to
$
y_{ref}(s)\thickapprox k_g c^{\top} H_o(s) r(s) \thickapprox D(s)
r(s).
$
We note that $k_g c^{\top} H_o(s)$ is exactly the reference model
of the MRAC design. Therefore this approach is equivalent to
mimicking MRAC, and, hence, high-gain feedback can be completely
avoided.

\begin{figure}[!h]
\begin{center}
\includegraphics[width=5.4in,height=1.2in]{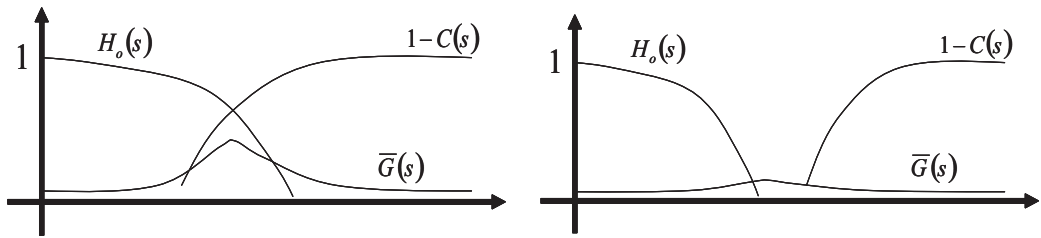}
\caption{ Cascaded systems. }\label{fig:nopass}
\end{center}
\end{figure}

However, increasing the bandwidth of $C(s)$ is not the only choice
for minimizing $ \Vert \bar{G}(s) \Vert_{\mathcal{L}_1}$.
Since $C(s)$ is a low-pass filter, its complementary $1-C(s)$ is a
high-pass filter with its cutoff frequency approximating the
bandwidth of $C(s)$. Since both $H_o(s)$ and $C(s)$ are strictly
proper systems, $\bar{G}(s)= H_o(s) (C(s)-1)$ is equivalent to
cascading a low-pass system $H_o(s)$ with a high-pass system
$C(s)-1$. If one chooses the cut-off frequency of $C(s)-1$ larger
than the bandwidth of $H_o(s)$, it ensures that $\bar{G}(s)$ is a
``no-pass'' system, and hence its   $\mathcal{L}_1$ gain can be
rendered arbitrarily small. This can be achieved via higher order
filter design methods. The illustration is given in Fig.
\ref{fig:nopass}.

To minimize $\Vert h_3 \Vert_{\mathcal{L}_\infty}$, we note that
$\Vert h_3 \Vert_{\mathcal{L}_\infty}$ can be upperbounded in two
ways:
$$
(i)~~~ \Vert h_3 \Vert_{\mathcal{L}_\infty} \leq \Vert (C(s)-1)
r(s)\Vert_{\mathcal{L}_1} \Vert h_4 \Vert_{\mathcal{L}_\infty}\,,
$$
where $h_4(t)$ is the  inverse Laplace transformation of $H_4(s)=
C(s) k_g H_o(s) \theta^{\top} H_o(s)$, and
$$ 
(ii) ~~~\Vert h_3 \Vert_{\mathcal{L}_\infty} \leq \Vert (C(s)-1)
C(s)\Vert_{\mathcal{L}_1} \Vert h_5 \Vert_{\mathcal{L}_\infty}\,,
$$ 
where $h_5(t)$ is the inverse Laplace transformation of $H_5(s)=
r(s) k_g H_o(s) \theta^{\top} H_o(s)$.

We note that since $r(t)$ is a bounded signal and $C(s), H_o(s)$
are stable  proper systems, $\Vert h_4 \Vert_{\mathcal{L}_\infty}$
and $\Vert h_5 \Vert_{\mathcal{L}_\infty}$ are finite.
 Therefore, $\Vert h_3 \Vert_{\mathcal{L}_\infty}$ can be minimized
 by
minimizing $\Vert (C(s)-1) r(s)\Vert_{\mathcal{L}_1}$ or $\Vert
(C(s)-1) C(s)\Vert_{\mathcal{L}_1}$. Following the same arguments
as above and assuming that $r(t)$ is  in low-frequency range, one
can choose the cut-off frequency of $C(s)-1$ to be larger than the
bandwidth of the reference signal $r(t)$ to minimize $\Vert
(C(s)-1) r(s)\Vert_{\mathcal{L}_1}$. For minimization of $ \Vert
C(s) (C(s)-1)\Vert_{\mathcal{L}_1}$ notice that if $C(s)$ is an
ideal low-pass filter, then $C(s) (C(s)-1)=0$ and hence $\Vert h_3
\Vert_{\mathcal{L}_\infty}=0$. Since an ideal low-pass filter is
not physically implementable, one can  minimize $ \Vert C(s)
(C(s)-1)\Vert_{\mathcal{L}_1}$ via appropriate  choice of $C(s)$.

The above presented approaches ensure that  $C(s)\approx 1$ in the
bandwidth of $r(s)$ and $H_o(s)$. Therefore it follows from
(\ref{ydes_def}) that
$ 
y_{des}(s) = C(s) k_g c^{\top} H_o(s) r(s) \thickapprox k_g
c^{\top} H_o(s) r(s),
$ 
which along with (\ref{L1_app2}) yields $ y_{des}(s) \thickapprox
D(s) r(s).
$

\begin{rem}
 From Corollary \ref{cor:1} and Lemma
\ref{lem:design} it follows that the ${\mathcal L}_1$ adaptive
controller can generate a system response   to track
(\ref{ydes_def}) and (\ref{udes_def}) both in transient and steady
state if we set the adaptive gain large and minimize $\lambda$ or
$\Vert h_3\Vert_{{\mathcal L}_\infty}$. Notice that $u_{des}(t)$
in  (\ref{udes_def}) depends upon the unknown parameter $\theta$,
while $y_{des}(t)$ in (\ref{ydes_def}) does not. This implies that
for different values of $\theta$, the
 $\mathcal{L}_1$ adaptive controller
 will generate different control signals (dependent on $\theta$)
 to ensure uniform system response (independent of $\theta$). This is natural, since different unknown
parameters imply different systems, and to have similar response
for different systems the control signals have to be different.
Here is the obvious advantage of the ${\mathcal L}_1$ adaptive
controller in a sense that it controls a partially known system as
an LTI feedback controller would have done if  the unknown
parameters were known. Finally, we note that if the term $k_g C(s)
C(s) \theta^{\top} H_o(s)$ is dominated by $k_g C(s) K^{\top}
H_o(s)$, then the
 controller in (\ref{udes_def}) turns into a robust one, and
 consequently the
${\mathcal L}_1$ adaptive controller degenerates into robust
design.
\end{rem}

\begin{rem}
It follows from (\ref{cor1_00}) that in the presence of large
adaptive gain the
 ${\mathcal L}_1$ adaptive controller and the closed-loop system
 state with it approximate $ u_{ref}(t), y_{ref}(t)$. Therefore, we conclude from
(\ref{ref_x})  that  $y(t)$ approximates  the response of the LTI
system $ c^{\top} (I-\bar{G}(s)\theta^{\top})^{-1} G(s) $ to the
input $r(t)$, hence its transient performance specifications, such
as overshoot and settling time, can be derived for every value of
$\theta$. If we further  minimize $\lambda$ or $\Vert
h_3\Vert_{{\mathcal L}_\infty}$, it follows from Lemma
\ref{lem:design} that $y(t)$ approximates  the response of the LTI
system $ C(s) c^{\top} H_o(s). $ In this case, the ${\mathcal
L}_1$ adaptive controller leads to uniform transient performance
of $y(t)$ independent of the value of the unknown parameter
$\theta$. For the resulting ${\mathcal L}_1$ adaptive control
signal one can characterize the transient specifications such as
its amplitude and rate change  for every $\theta\in \Omega$, using
$u_{des}(t)$ for it.
\end{rem}

\section{Discussion}\label{sec:comp}

\subsection{Comparison to high-gain controller}

We use  a scalar system to compare the performance of the
$\mathcal{L}_1$ adaptive controller and a linear high-gain
controller. Towards that end, let $ \dot{x}(t) = -\theta x(t) +
u(t)\,, $ where $x\in \IR$ is the measurable system state, $u\in
\IR$ is the control signal and $\theta\in \IR$ is unknown, which
belongs to a given compact set $[\theta_{\min},\;\theta_{\max}]$.
 Let $u(t)=-k x(t)+k r(t)$, leading to the following
closed-loop system $ \dot{x}(t) = (-\theta-k) x(t)+ k r(t). $ We
need to choose $ k>-\theta_{\min}$ to guarantee  stability. We
note that both the steady state error and the transient
performance depend on the unknown parameter value $\theta$. By
further introducing a proportional-integral controller, one can
achieve zero steady state error. If one chooses $k\gg
\max\{|\theta_{\max}|,|\theta_{\min}| \}$, it leads to high-gain
system
$$ 
x(s) = \frac{k}{s-(-\theta-k)} r(s) \approx \frac{k}{s+k} r(s)\,.
$$ 

To apply the ${\mathcal L}_1$ adaptive controller, let the desired
reference system be
$ 
D(s)=\frac{2}{s+2}.
$ 
 Let
$ u_1= - 2 x, k_g=2, $ leading to $ H_o(s) = \frac{1}{s+ 2}. $
Choose $C(s)$ as in (\ref{firstorderCs}) with large $\omega_n$,
and set adaptive gain $\Gamma_c$ large. Then it follows from
Theorem \ref{thm:5} that
\begin{eqnarray}
x(s) & \approx & x_{ref}(s)=C(s) k_g H_o(s) r(s) \approx
\frac{\omega_n}{s+\omega_n} \frac{2}{s+2} r(s) \approx \frac{2}{s+2} r(s)\label{1} \\
u(s) &  \approx & u_{ref}(s) = (- 2 +\theta) x_{ref}(s) + 2
r(s).\label{2}
\end{eqnarray}
The relationship in (\ref{1}) implies that the control objective
is met, while the relationship in (\ref{2}) states that the
${\mathcal L}_1$ adaptive controller approximates $u_{ref}(t)$,
which cancels the unknown $\theta$.

\subsection{Time-delay margin in the presence of large adaptive gain}

A well-known fact in robust control is that the high gain in the
feedback loop can lead to increased control effort and reduced
phase margin. Since we argue that the performance bounds of
$\mathcal L_1$ adaptive controller can be systematically improved
by increasing the adaptation gain, in this section we provide a
brief robustness analysis of the $\mathcal L_1$ adaptive
controller in parallel to MRAC. To enable the use of the frequency
domain tools for robustness analysis, we consider a  scalar linear
system in the presence of constant unknown disturbance and close
the loop with a MRAC controller and $\mathcal L_1$ controller. So,
let $ \dot{x}(t) = x(t)+ u(t)+\theta$,
 where $x\in
\IR$ is the measured  state, $u\in \IR$ is the control signal,
$\theta\in \IR$ is an unknown constant parameter. If we apply the
MRAC controller, then it reduces to the well-known PI structure:
$$ u(t) = -\hat{\theta}(t)-2 x(t)+r(t)\,,\quad
\dot{x}_m(t)=-x_m(t)+r(t)\,,\quad \dot{\hat{\theta}}(t) =\Gamma
(x(t)-x_m(t))\,. $$ The open-loop transfer function for the
time-delay margin analysis of this controller in the presence of
the time-delay at the system input is $
H_o(s)=\frac{-ks+\Gamma}{s(s-1)} e^{-s \tau} $. Application of the
$\mathcal L_1$ controller leads to a filtered version of the PI
controller:
\begin{eqnarray}
u(s)&=&- C(s) (\hat{\theta}(s)-r(s))-2 x(s), \qquad
\dot{\hat{\theta}}(t) =-\Gamma (\hat{x}(t)-x(t))\,\nonumber \\
\dot{\hat{x}}(t) & =& -\hat{x}(t) +u(t)+\hat{\theta}(t) .\nonumber
\end{eqnarray}
In this case, the open-loop transfer function for the time-delay
margin analysis
 in the presence of the time-delay at the system input is $H_o(s)
= \frac{-k}{s-1} e^{-s \tau}+\frac{\Gamma C(s)}{s^2-a_m s +
\Gamma} (e^{-s \tau}-1) $.
\begin{figure}[!h]
\begin{center}
\includegraphics[width=2.8in,height=1.8in]{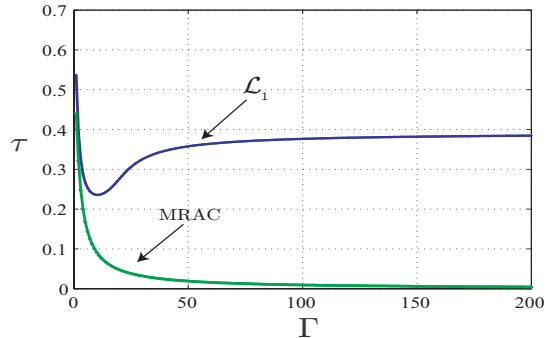}
\caption{Effects of adaptive gain on time-delay margin in MRAC and
$\mathcal L_1$ adaptive controller with $C(s)=\frac{s}{s+1}$
}\label{fig:tdmargin}
\end{center}
\end{figure}
We plot the time-delay margin of both systems with respect to
 adaptive gain $\Gamma$ in Fig. \ref{fig:tdmargin}. We notice that the time-delay
margin of PI controller goes to zero as $\Gamma\rightarrow
\infty$, while the time-delay margin of $\mathcal L_1$ adaptive
control architecture is bounded away from zero as
$\Gamma\rightarrow \infty$. Details on this analysis can be found
in \cite{chengyu_gnc06}.



\section{Time-varying unknown parameters}\label{sec:timevarying}

 In this section, we consider the performance of the
$\mathcal{L}_1$ adaptive controller in the presence of
time-varying unknown parameters. We prove that, in this case as
well, by increasing the adaptation gain one can ensure uniform
transient response for system's both signals, input and output,
simultaneously. We consider the same system in (\ref{problem})
with unknown time-varying parameters $\theta(t)\in \IR^n$,
assuming that
$
\theta(t) \in \Theta,\; \forall\; t\geq 0\,.
$
 We further
assume that $\theta(t)$ is continuously differentiable with
uniformly bounded derivative:
\begin{equation}
\Vert \dot{\theta}(t) \Vert  \leq  d_{\theta}<\infty, \quad
\forall~ t\geq 0\,
 ,\label{derivativeB}
\end{equation}
where the number $d_{\theta}$ can be arbitrarily large. We
consider the same reference system in (\ref{ref_u}) with $\eta(t)$
defined as
\begin{equation}\label{eta_newdef}
\eta(t) = \theta^{\top}(t) x_{ref}(t)\,.
\end{equation}
Hence, (\ref{yrefdef}) becomes
\begin{equation}\label{yrefdefnew}
x_{ref}(s)  =  H_o(s) \left(k_g C(s) r(s)+ (C(s)-1)
\eta(s)\right)\,,
\end{equation}
where $\eta(s)$ is the Laplace transformation of $\eta(t)$ in
(\ref{eta_newdef}). Let $\eta_1(t)$ be the signal with its Laplace
transformation given by
\begin{equation}
\eta_1(s) = H_o(s) (C(s)-1) \eta(s)\,.
\end{equation}
It can be derived easily that
\begin{equation}
\Vert \eta_1 \Vert_{\mathcal{L}_\infty} \leq \Vert H_o(s)
(C(s)-1)\Vert_{\mathcal{L}_1} \theta_{\max}\Vert x_{ref}
\Vert_{\mathcal{L}_\infty} \,,
\end{equation}
where $\theta_{\max}$ is defined in (\ref{thetabardef}). It
follows from  Theorem \ref{thm:sg} that the closed-loop reference
system is stable if the same $\mathcal{L}_1$-gain requirement in
(\ref{condition3}) holds. Instead of Lemma \ref{lem:1} and Theorem
\ref{thm:5}, we have the following results.
\begin{lem}\label{lem:newbounds}
For the system in (\ref{problem}) in the presence of unknown
time-varying $\theta(t)$, we have
\begin{equation}
\Vert \tilde{x} \Vert_{\mathcal{L}_{\infty}}  \leq
\sqrt{\theta_{m}/(\lambda_{\min}(P) \Gamma_c)}\,,
\label{barthemax_time}
\end{equation}
where
\begin{eqnarray}
 \quad \theta_{m}  \triangleq  \max_{\theta\in \Theta}
\sum_{i=1}^{n} 4 \theta_i^2 +2\frac{d_{\theta} \lambda_{\max}(P)
}{\lambda_{\min}(Q)}
  \max_{\theta\in \Theta}  \Vert\theta\Vert
 \,.& &\label{thetamdef}
\end{eqnarray}
\end{lem}
{\bf Proof.}  Using the same candidate Lyapunov function in
(\ref{lyapunov}), it follows that
\begin{equation}\label{lem1_11}
\dot{V}(t) \leq -\tilde{x}^{\top}(t) Q \tilde{x}(t) +2
\Gamma_{c}^{-1} \tilde{\theta}^{\top}(t) \dot{\theta}(t)\,.
\end{equation}
If at any $t$,
\begin{equation}\label{lem1new_15}
V(t)>\theta_{m}/\Gamma_c\,,
\end{equation}
where $\theta_{m}$ is defined in (\ref{thetamdef}), then it
follows from (\ref{lem1_5}) that
\begin{equation}
\tilde{x}^{\top}(t) P \tilde{x}(t) >2 \frac{d_{\theta}
\lambda_{\max}(P)}{\Gamma_c \lambda_{\min}(Q)}
  \max_{\theta\in \Theta}  \Vert\theta\Vert  ,
\end{equation}
and hence
\begin{eqnarray}
\tilde{x}^{\top}(t) Q \tilde{x}(t)  >  \frac{\lambda_{\min}(Q)}{
\lambda_{\max}(P)} \tilde{x}^{\top}(t) P \tilde{x}(t)
 > 2 \frac{
  d_{\theta}\displaystyle{ \max_{\theta\in \Theta} } \Vert\theta\Vert
 }{\Gamma_c}\,.& & \nonumber
\end{eqnarray}
The upper bounds in  (\ref{derivativeB}) along with the projection
based adaptive laws lead to the following upper bound:
$$
 ( 2\tilde{\theta}^{\top}(t)
\dot{\theta}(t))/\Gamma_c\leq  \frac{2
d_{\theta}\displaystyle{\max_{\theta\in \Theta} \Vert\theta\Vert}
  }{\Gamma_c}\,.
$$
Hence,  it follows from (\ref{lem1_11}) and (\ref{lem1new_15})
that
\begin{equation}\label{lem1_21}
\dot{V}(t) <0\,.
\end{equation}
 Since $
V(0) \leq \theta_{m}/\Gamma_c$, it follows from (\ref{lem1_21})
that $ V(t) \leq \theta_{m}/\Gamma_c$ for any $t \geq 0$. Since $
\lambda_{\min}(P) \Vert \tilde{x}(t) \Vert^2 \leq
\tilde{x}^{\top}(t) P \tilde{x}(t)\leq V(t)$, then
$$ 
|| \tilde{x}(t) ||^2  \leq \frac{\theta_{m}}{\lambda_{\min}(P)
\Gamma_c}\,,
$$ 
which concludes the proof.
 $\hfill{\square}$

\begin{thm}\label{thm:8}
Given the system in (\ref{problem}) with unknown time-varying
$\theta(t)$ and the  $\mathcal{L}_1$ adaptive controller defined
via (\ref{u1_design}), (\ref{L1_companionmodal}),
(\ref{adaptivelaw_L1}), (\ref{us_def}) subject to
(\ref{condition3}), we have:
\begin{eqnarray}
\Vert x-x_{ref} \Vert_{{\mathcal L}_{\infty}}  & \leq & \gamma_3\,
,\label{thm8_00} \\
\Vert u - u_{ref} \Vert_{{\mathcal L}_{\infty}} & \leq &
\gamma_4\, ,\label{thm8_01}
\end{eqnarray}
where
\begin{eqnarray}
\gamma_3 & = & \frac{\Vert C(s) \Vert_{\mathcal{L}_1}}{1-\Vert
H_o(s) (1-C(s))\Vert_{\mathcal{L}_1} \theta_{\max}}
\sqrt{\frac{\theta_{m}}{\lambda_{\max}(P) \Gamma_c}}\,,
\label{gamma3def} \\
\gamma_4 & = & \Big\|C(s) \frac{1}{c_{o}^{\top} H_o(s)}
c_{o}^{\top}
\Big\|_{\mathcal{L}_1}\sqrt{\frac{{\theta}_{m}}{\lambda_{\max}(P)\Gamma_c}}+\left(\Vert
K^{\top}\Vert +\Vert C(s) \Vert_{\mathcal{L}_1} \theta_{\max}
\right) \gamma_3\,. \label{gamma4def}
\end{eqnarray}
\end{thm}
\medskip

{\bf Proof.}
Letting $ \tilde{r}(t)  = \tilde{\theta}^{\top}(t) x(t) \,, $
$\eta_2(t)  =  \theta^{\top}(t) x(t) $, it follows from the system
in (\ref{problem}) and the control law in
(\ref{u1_design})-(\ref{us_def}) that
\begin{equation}\label{thm1_5}
x(s) = H_o(s)  \left( (C(s)-1) \eta_2(s) +C(s) k_g r(s)+ C(s)
\tilde{r}(s)\right)\,.
\end{equation}
Following the
definition of $r_2(t)$ in (\ref{thm3_3new}), it follows from
(\ref{yrefdefnew}) and (\ref{thm1_5}) that
\begin{equation}
r_2(s) = H_o(s) \left( (C(s)-1) r_3(s) -   C(s)
\tilde{r}(s)\right), \quad r_2(0)=0\,,\label{thm1_7}
\end{equation}
where $r_3(s)$ is the Laplace transformation of the signal
\begin{equation}\label{thm1_11}
r_3(t) = \theta^{\top}(t) r_2(t)\,.
\end{equation}
Lemma \ref{lem:L1} gives the following upper bound:
\begin{equation}\label{thm1_8}
\Vert r_{2_t} \Vert_{\mathcal{L}_{\infty}} \leq \Vert H_o(s)
(1-C(s))\Vert_{\mathcal{L}_1}  \Vert r_{3_t}
\Vert_{\mathcal{L}_{\infty}} + \Vert r_{4_t}
\Vert_{\mathcal{L}_{\infty}}\,,
\end{equation}
where $ r_4(t)$ is the signal with its Laplace transformation
$
r_4(s) = C(s) H_o(s) \tilde{r}(s).
$
Since
$
\tilde{x}(s) = H_o(s) \tilde{r}(s)\,,
$
we have
$
r_4(s) = C(s) \tilde{x}(s)\,,
$
and hence
$
\Vert r_{4_t} \Vert_{\mathcal{L}_{\infty}} \leq \Vert C(s)
\Vert_{\mathcal{L}_1} \Vert \tilde{x}_t
\Vert_{\mathcal{L}_{\infty}}\,.
$
Using the definition of $\theta_{\max}$ in (\ref{thetabardef}),
one can easily  verify from (\ref{thm1_11})
 that
$
\Vert r_{3_t} \Vert_{\mathcal{L}_{\infty}} \leq \theta_{\max}
\Vert r_{2_t} \Vert_{\mathcal{L}_{\infty}}\,.
$
From (\ref{thm1_8}) we have
\begin{equation}
\Vert r_{2_t} \Vert_{\mathcal{L}_{\infty}} \leq \Vert H_o(s)
(1-C(s))\Vert_{\mathcal{L}_1} \theta_{\max} \Vert r_{2_t}
\Vert_{\mathcal{L}_{\infty}}+  \Vert C(s)
\Vert_{\mathcal{L}_1}\Vert \tilde{x}_t
\Vert_{\mathcal{L}_{\infty}}\,.
\end{equation}
The upper bound from Lemma \ref{lem:newbounds} and the $\mathcal
L_1$-gain requirement from (\ref{condition3}) lead to the
following upper bound
\begin{equation}
\Vert r_{2_t} \Vert_{\mathcal{L}_{\infty}} \leq \frac{\Vert C(s)
\Vert_{\mathcal{L}_1}}{1-\Vert H_o(s)
(1-C(s))\Vert_{\mathcal{L}_1} \theta_{\max} }
\sqrt{\frac{\theta_{m}}{\lambda_{\max}(P) \Gamma_c}}\,,
\end{equation}
which holds uniformly for all $t\ge 0$ and therefore leads to
(\ref{thm8_00}).

To prove the bound in (\ref{thm8_01}), we notice that from
(\ref{u1_design})-(\ref{us_def}) and (\ref{ref_u})-(\ref{eta_def})
one can derive
\begin{equation}\label{lem2_new5}
u(s)-u_{ref}(s) = - K^{\top}(x(s)-x_{ref}(s)) +C(s)
(\eta_2(s)-\eta(s)) +r_5(s)\,,
\end{equation}
where $r_5(s)= C(s)\tilde{r}(s)$. Therefore, it follows from Lemma
\ref{lem:L1} that
\begin{equation}\label{thm1_99}
\Vert u-u_{ref}\Vert_{\mathcal{L}_{\infty}} \leq \left(\Vert
K^{\top}\Vert +\Vert C(s) \Vert_{\mathcal{L}_1} \theta_{\max}
\right)\Vert x-x_{ref}\Vert_{\mathcal{L}_{\infty}}+\Vert
r_5\Vert_{\mathcal{L}_{\infty}}\,.
\end{equation}
We have $ r_5(s)  = C(s) \frac{1}{c_{o}^{\top} H_o(s)}
c_{o}^{\top} H_o(s) \tilde r(s)       =  C(s)
\frac{1}{c_{o}^{\top} H_o(s)} c_{o}^{\top} \tilde{x}(s)\,, $ and
hence,
$$ 
\Vert r_5\Vert_{{\mathcal L}_{\infty}} \leq \Big\| C(s)
\frac{1}{c_{o}^{\top} H(s)} c_{o}^{\top} \Big\|_{\mathcal{L}_1}
\Vert \tilde{x} \Vert_{{\mathcal L}_{\infty}}\,.
$$ 
 Lemma \ref{lem:newbounds}  consequently leads to the upper bound:
$$ 
\Vert r_5\Vert_{{\mathcal L}_{\infty}} \leq  \Big\| C(s)
\frac{1}{c_{o}^{\top} H(s)} c_{o}^{\top} \Big\|_{\mathcal{L}_1}
\sqrt{\frac{\theta_{m}}{\lambda_{\max}(P) \Gamma_c}}\,,
$$ 
which, when substituted into (\ref{thm1_99}), leads to
(\ref{thm8_01}).
 $\hfill{\square}$

Since (\ref{condition3}) ensures the stability of the reference
system, it follows from Theorem \ref{thm:8} that the same
$\mathcal{L}_1$-gain requirement ensures the stability of
$\mathcal{L}_1$ adaptive controller. Theorem \ref{thm:8} further
implies that the $\mathcal L_1$ adaptive controller approximates
$u_{ref}(t)$ both in transient and steady state. It is
straightforward to verify that Corollary \ref{cor:1} holds for
time-varying unknown $\theta(t)$ as well.

 We note that the
control law $u_{ref}(t)$ in the closed-loop reference system,
which is used in the analysis of $\mathcal L_\infty$ norm bounds,
is not implementable since its definition involves the unknown
parameters.   So, it is important to understand how these bounds
can be used for ensuring uniform transient response with {\em
desired} specifications. We notice that the following {\em ideal}
control signal
\begin{equation}\label{uideal}
u_{ideal}(t)=k_g r(t)+\theta^{\top}(t) x_{ref}(t) -K^{\top}
x_{ref}(t)
\end{equation}
is the one that leads  to desired system response:
\begin{eqnarray}
\dot{x}_{ref}(t) &=& A_m x_{ref}(t)+b k_g r(t)\label{desrefstate}\\
y_{ref}(t)&=&c^\top x_{ref}(t)\label{desrefoutput}
\end{eqnarray}
by cancelling the uncertainties exactly. If a part of
 $u_{ideal}(t)$ is  low-pass filtered by $C(s)$  in
(\ref{ref_u}), then $u_{ref}(t)$  cancels  the uncertainties
dependent upon the bandwidth of $C(s)$. In case of fast varying
$\theta(t)$, it is obvious that the bandwidth of the controller
needs to be matched correspondingly.

\section{Simulations}\label{sec:simu}
Consider the system in (\ref{problem}) with the following
parameters:
$$ A=\left[
\begin{array}{lr}
0 & 1 \\
-1 & -1.4
\end{array}
\right]\,, \quad b=[0 \quad 1]^{\top} ,\quad  c=[1 \quad
0]^{\top}\,, \quad \theta=[4\quad -4.5]^{\top}\,. $$ We further
assume that the unknown parameter $\theta$ belongs to a known
compact set $ \Theta=\{\theta\in \IR^2 \;|\; \theta_1\in[-10,10],
\theta_2\in[-10,10]\}$.

\begin{figure}[h!]
\centering \mbox{\subfigure[ $\lambda$ (solid) defined in
(\ref{simu_omega})
 ]{\epsfig{file=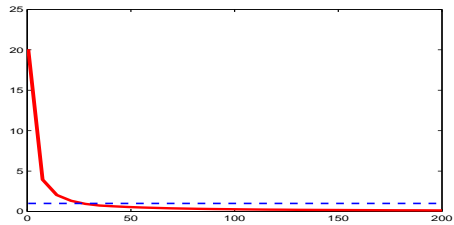,width=2.8in,height=1.3in
}\label{fig:criterion}}\quad \subfigure[ $\lambda$ (solid) defined
in (\ref{simu_omega2}) ]{\epsfig{file=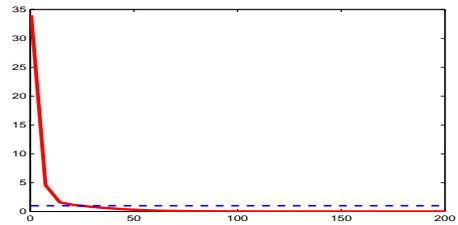,width=2.8in,
height=1.3in}\label{fig:criterion2} }} \caption{$\lambda$ (solid)
 with respect to $\omega$ and
constant $1$ (dashed) }
\end{figure}

We give now the complete $\mathcal{L}_1$ adaptive controller for
this system. Since $A$ is Hurwitz, we set $ K = 0$. Letting
$\Gamma_c = 10000 $, we implement the $L_1$ adaptive controller
following (\ref{u1_design}), (\ref{L1_companionmodal}),
(\ref{adaptivelaw_L1}) and (\ref{us_def}). First,   we check
stability of this $\mathcal{L}_1$ adaptive controller. It follows
from (\ref{thetabardef}) that $ \theta_{\max} =20 $ and $\Vert
\bar{G}\Vert_{L_1}$ can be calculated numerically. In Fig.
\ref{fig:criterion}, we plot
\begin{equation}\label{simu_omega}
\lambda=\Vert \bar{G}\Vert_{L_1} \theta_{\max}
\end{equation}
 with respect to
$\omega$ and compare it to $1$. We notice that for $\omega
> 30$, we have $\lambda<1$, and  the
$\mathcal{L}_1$ gain requirement for stability is satisfied. So,
we can choose
\begin{eqnarray}\label{simu:Cs}
C(s) & = & \frac{160}{s+160}
\end{eqnarray}
to ensure that $\lambda < 0.01$, which consequently leads to
improved performance bounds in
(\ref{lemdesign_00})-(\ref{lemdesign_02}). For $\omega=160$,
 we have
$ \lambda =\Vert \bar{G}(s) \Vert_{\mathcal{L}_1} \theta_{\max}  =
0.1725<1, $ so the $\mathcal{L}_1$-gain requirement  in
(\ref{condition3}) is indeed satisfied.

\begin{figure}[h!]
\centering \mbox{\subfigure[ $y(t)$ (solid) and $r(t)$ (dashed)
 ]{\epsfig{file=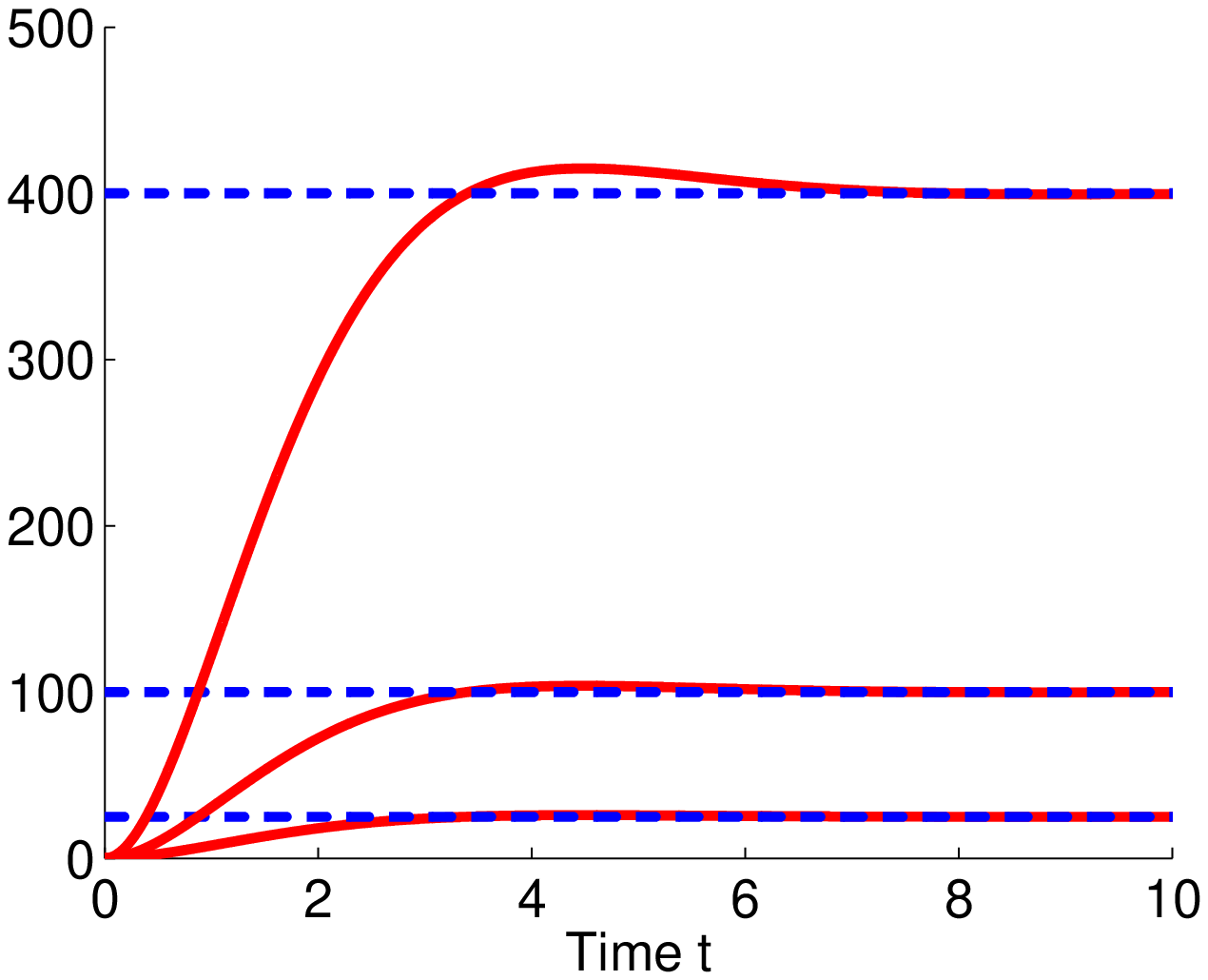,width=2.8in,height=1.3in
}\label{fig:L1_y}}\quad \subfigure[ Time-history of $u(t)$
]{\epsfig{file=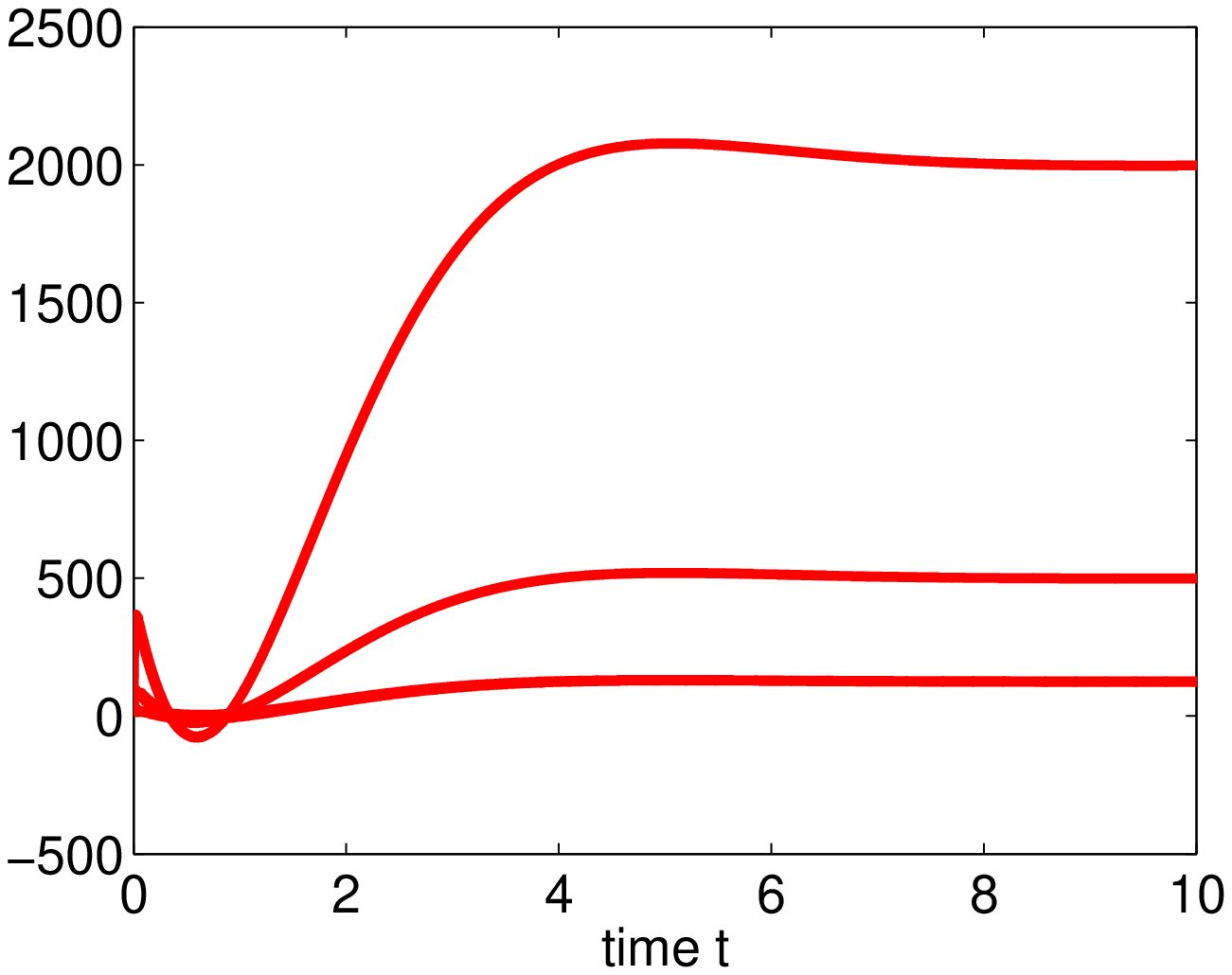,width=2.8in, height=1.3in}\label{fig:L1_u}
}} \caption{Performance of $\mathcal{L}_1$ adaptive controller
with $C(s)=\frac{160}{s+160}$ for $r=25,~100,~400 $  }
\end{figure}

\begin{figure}[h!]
\centering \mbox{\subfigure[ $y(t)$ (solid) and $r(t)$ (dashed)
  ]{\epsfig{file=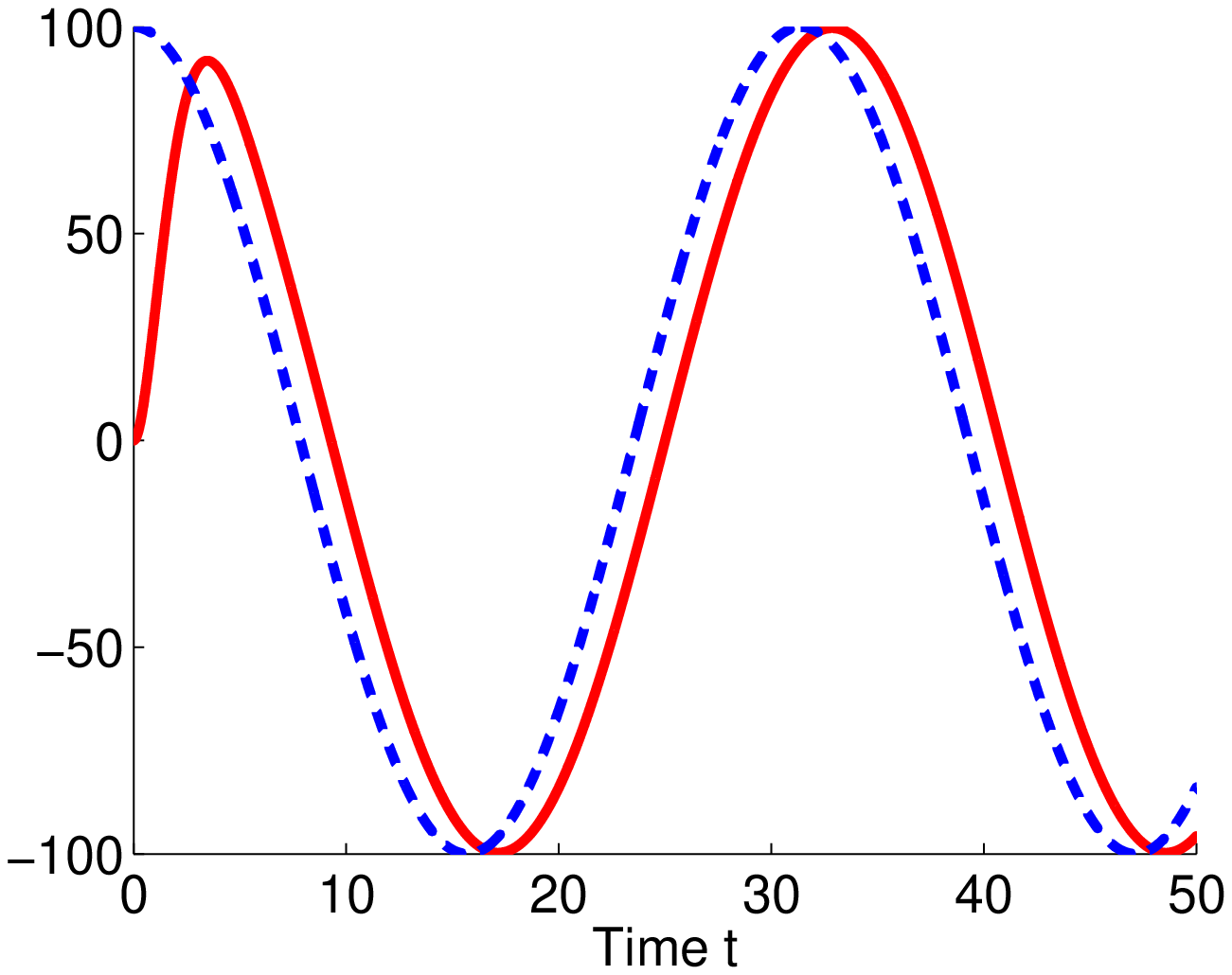,width=2.8in,height=1.3in
}\label{fig:siny1}}\quad \subfigure[ Time-history of $u(t)$
]{\epsfig{file=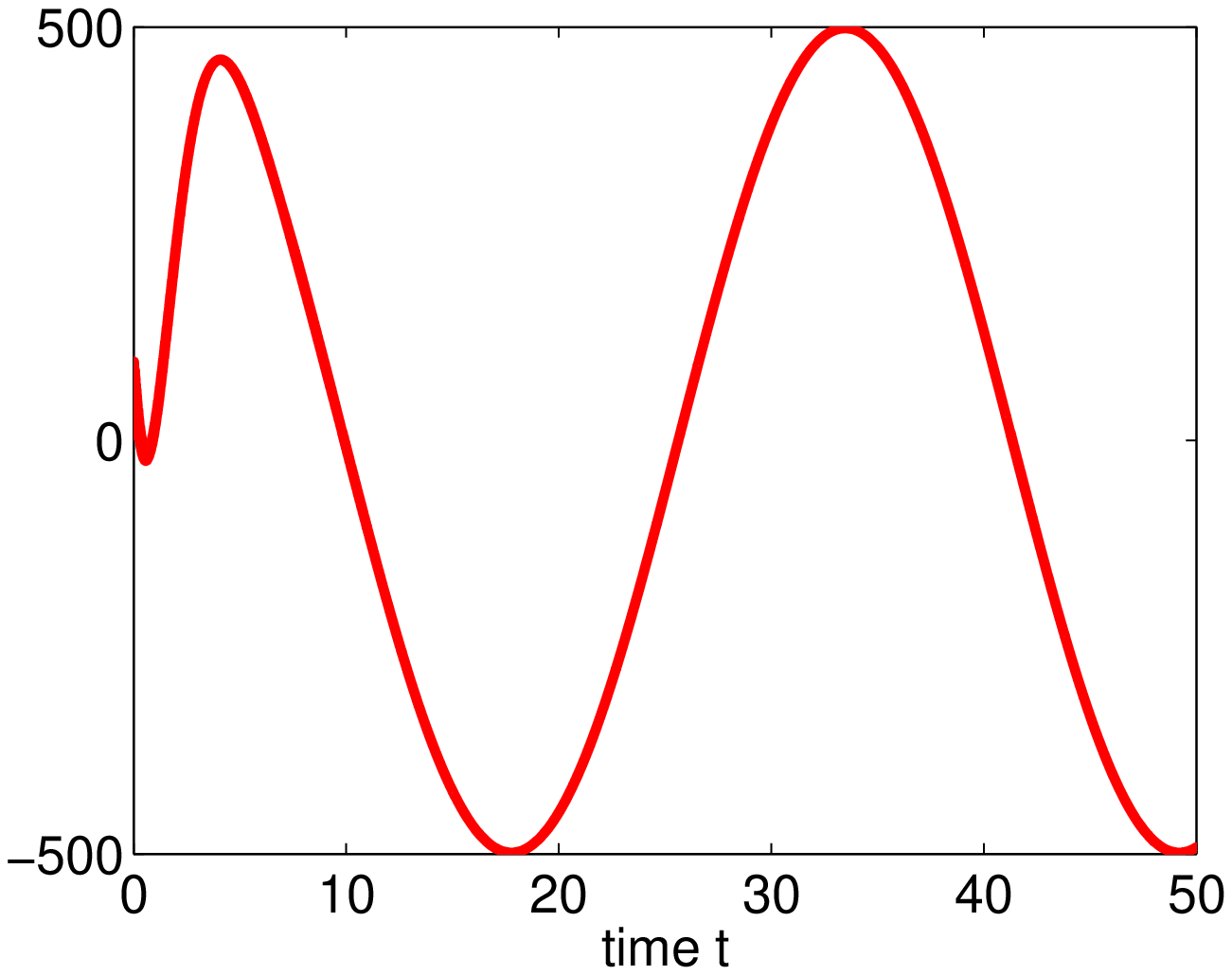,width=2.8in,
height=1.3in}\label{fig:sinu1}}} \caption{ Performance of
$\mathcal{L}_1$ adaptive controller with $C(s)=\frac{160}{s+160}$
for $r=100 \cos (0.2 t)$ }
\end{figure}

\begin{figure}[h!]
\centering \mbox{\subfigure[ $y(t)$ (solid) and $r(t)$ (dashed)
  ]{\epsfig{file=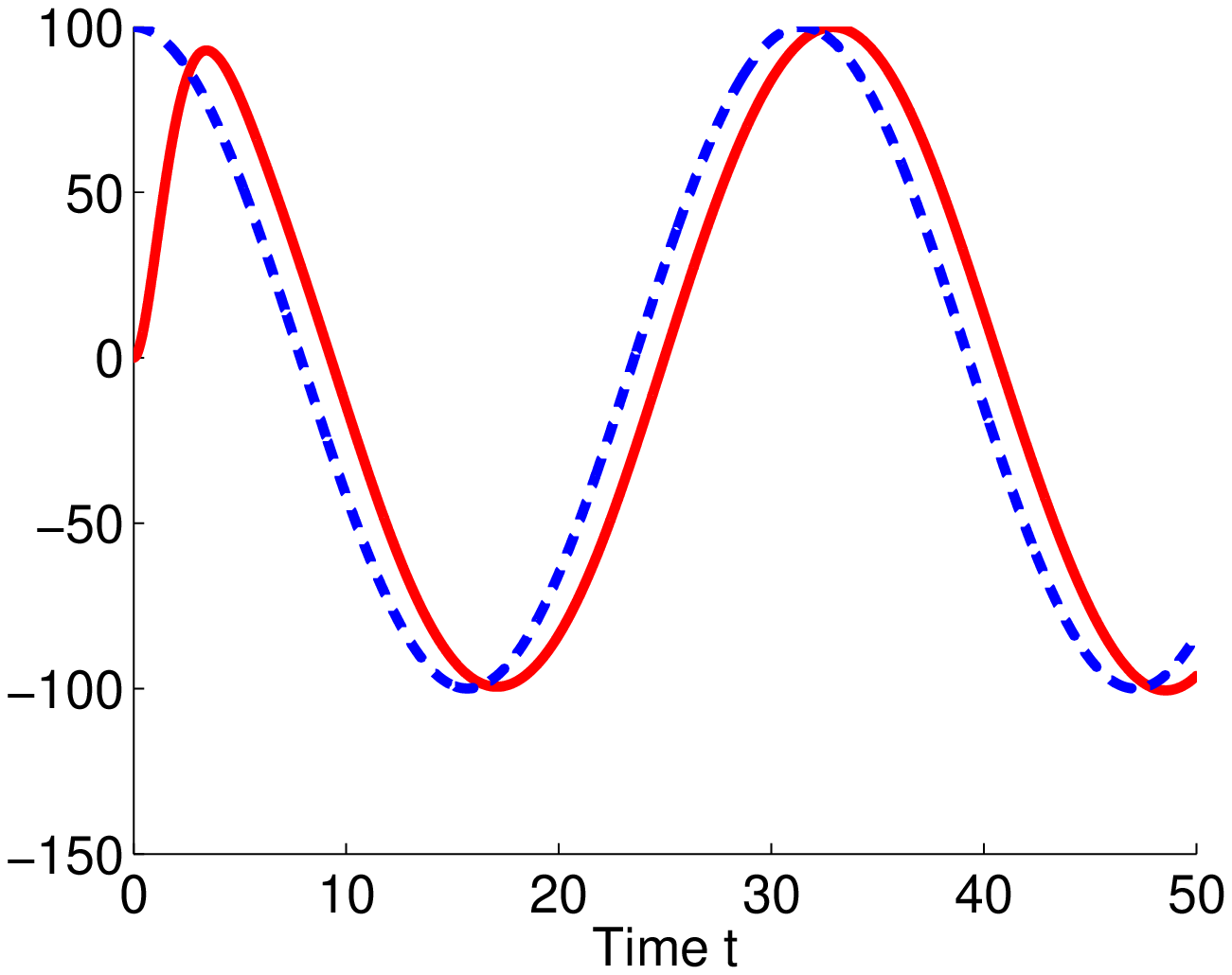,width=2.8in,height=1.3in
}\label{fig:siny1tv}}\quad \subfigure[ Time-history of $u(t)$
]{\epsfig{file=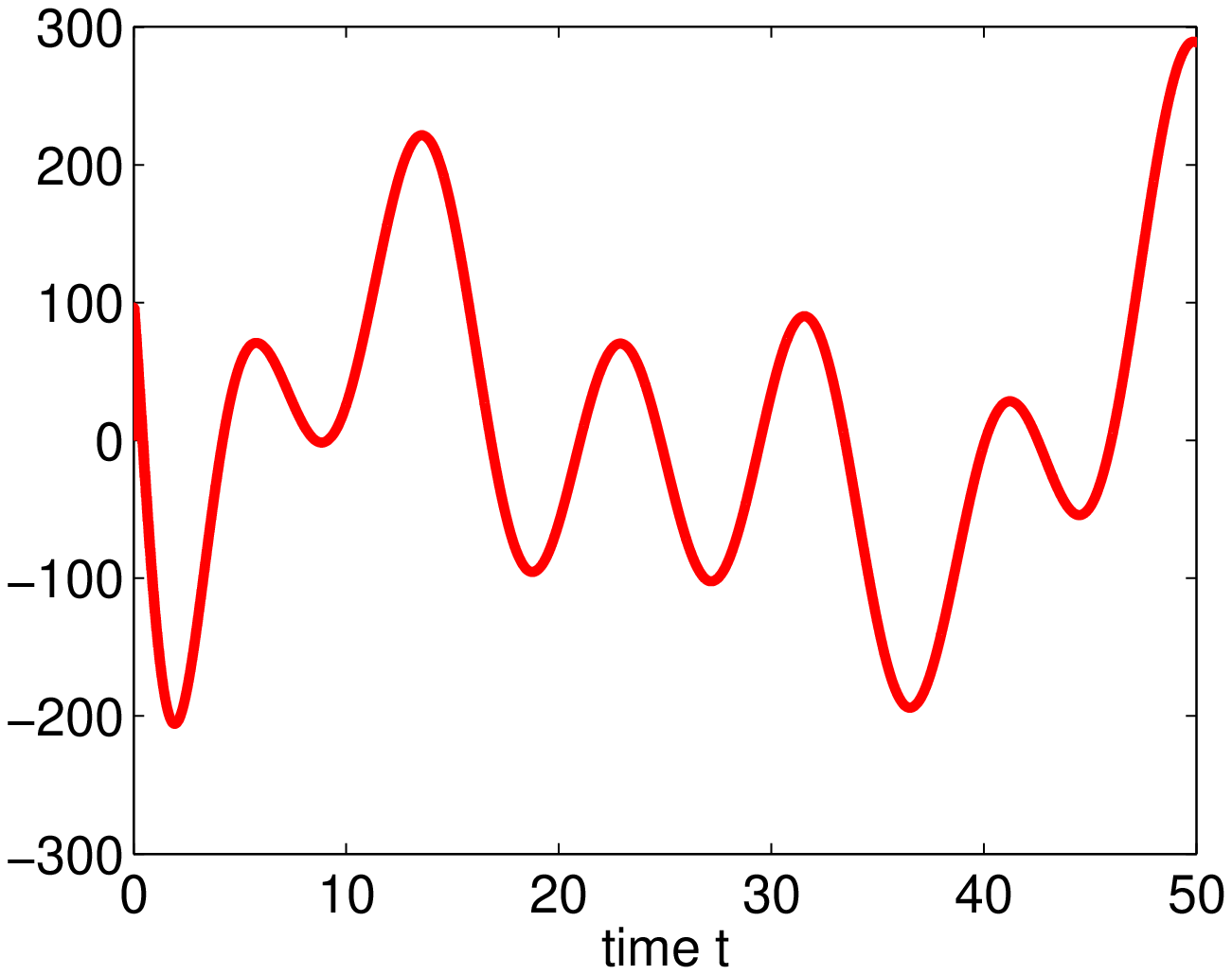,width=2.8in,
height=1.3in}\label{fig:sinu1tv}}} \caption{ Performance of
$\mathcal{L}_1$ adaptive controller with $C(s)=\frac{160}{s+160}$
for $r=100 \cos (0.2 t)$ with time-varying $\theta(t)=[2+2\cos(0.5
t)\quad 2+0.3\cos(0.5 t)+0.2\cos(t/\pi)]^{\top}$}
\end{figure}

 The simulation results
of the $\mathcal{L}_1$ adaptive controller are shown in Figs.
\ref{fig:L1_y}-\ref{fig:L1_u} for reference inputs
$r=25,~100,~400$, respectively. We note that it leads to scaled
control inputs and scaled system outputs for scaled reference
inputs.
Figs. \ref{fig:siny1}-\ref{fig:sinu1} show the system response and
the control signal for  reference input $r(t) = 100 \cos (0.2 t)$,
without any retuning of the controller.  Figs.
\ref{fig:siny1tv}-\ref{fig:sinu1tv} show the system response and
the control signal for  reference input $r(t) = 100 \cos (0.2 t)$
and time varying $\theta(t)=[2+2\cos(0.5 t)\quad 2+0.3\cos(0.5
t)+0.2\cos(t/\pi)]^{\top}$, without any retuning of the
controller. We note that the $\mathcal{L}_1$ adaptive controller
leads to almost identical tracking performance for both constant
or time-varying unknown parameters. The control signals are
different since they are adapting to different uncertainties to
ensure uniform transient response.


\begin{figure}[h!]
\centering \mbox{\subfigure[ $y(t)$ (solid) and $r(t)$ (dashed)
 ]{\epsfig{file=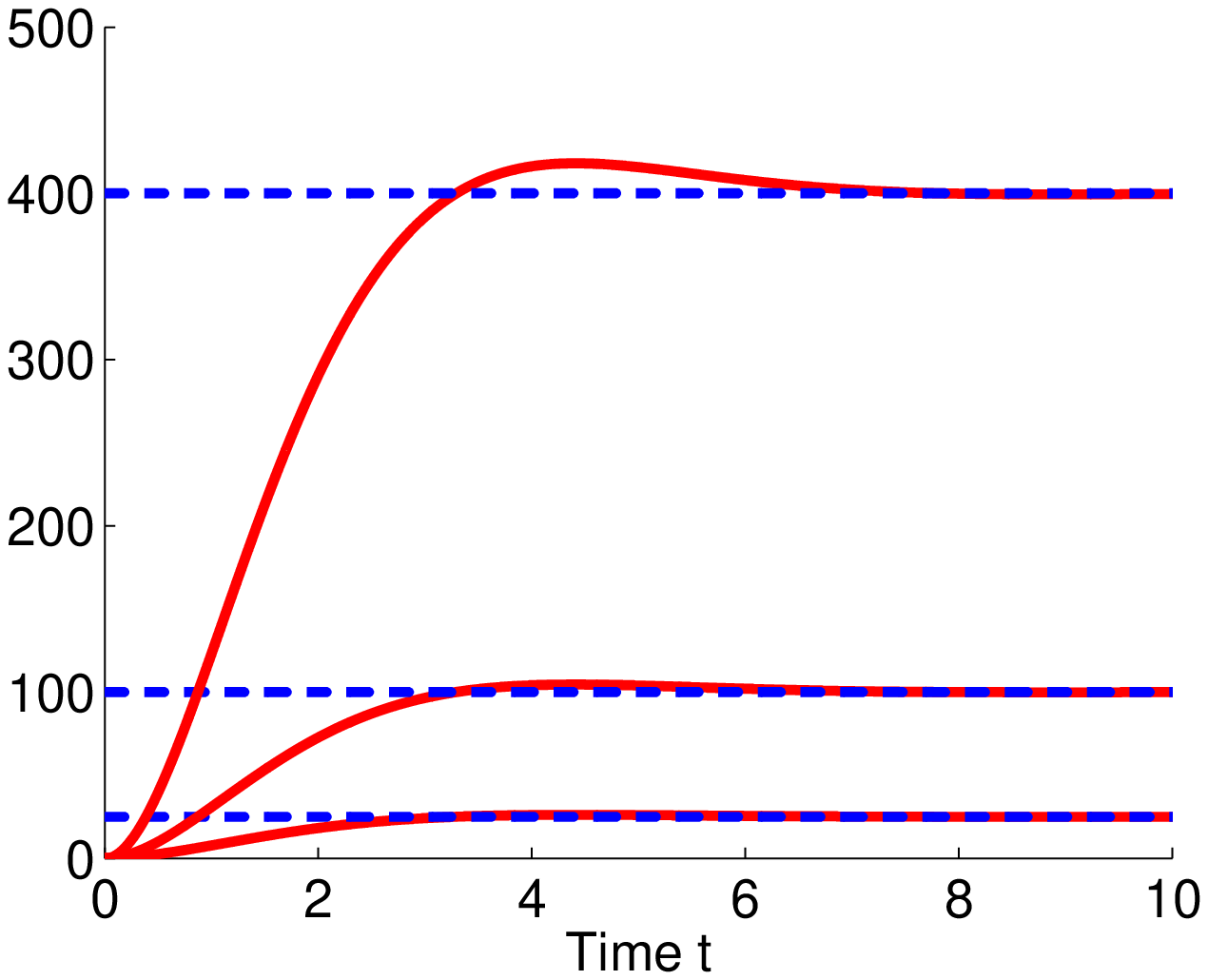,width=2.8in,height=1.3in
}\label{fig:L1_y_r200}}\quad \subfigure[ Time-history of $u(t)$
]{\epsfig{file=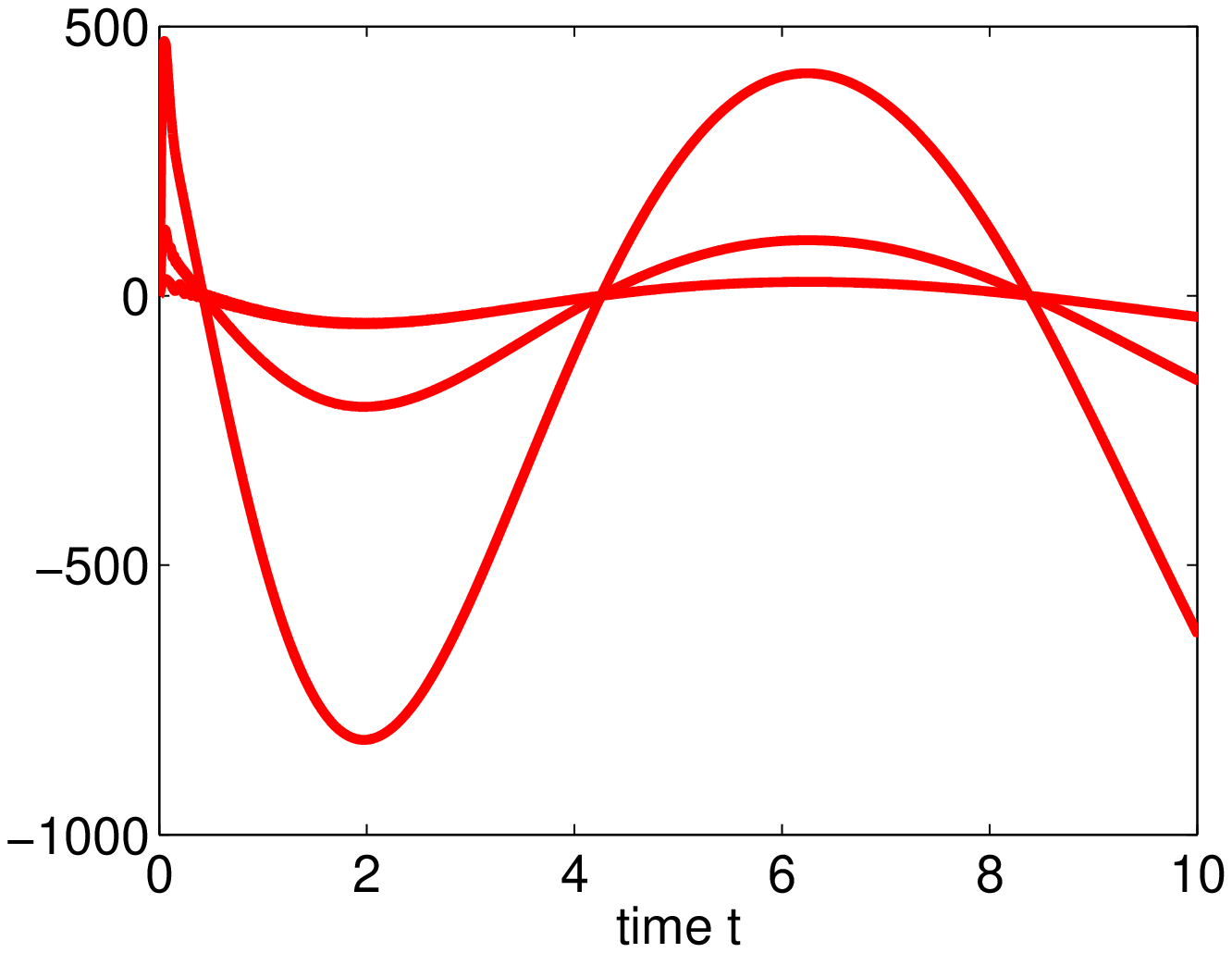,width=2.8in,
height=1.3in}\label{fig:L1_u_r200}}}\caption{Performance of
$\mathcal{L}_1$ adaptive controller with $C(s)=\frac{7500 s
+50^3}{(s+50)^3}$ for $r=25,~100,~400 $  }
\end{figure}

\begin{figure}[h!]
\centering \mbox{\subfigure[ $y(t)$ (solid) and $r(t)$ (dashed)
 ]{\epsfig{file=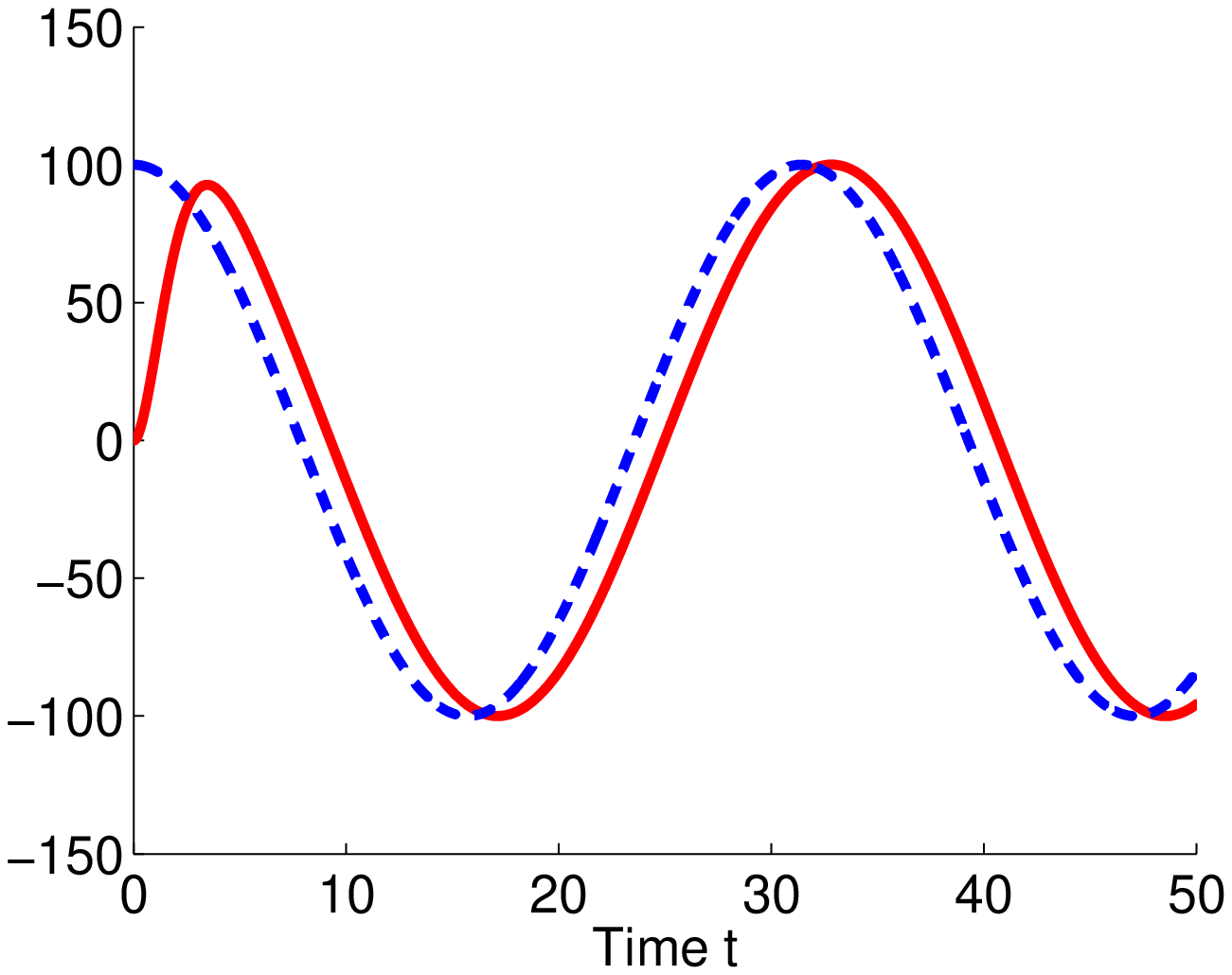,width=2.8in,height=1.3in
}\label{fig:siny2}}\quad \subfigure[ Time-history of $u(t)$
]{\epsfig{file=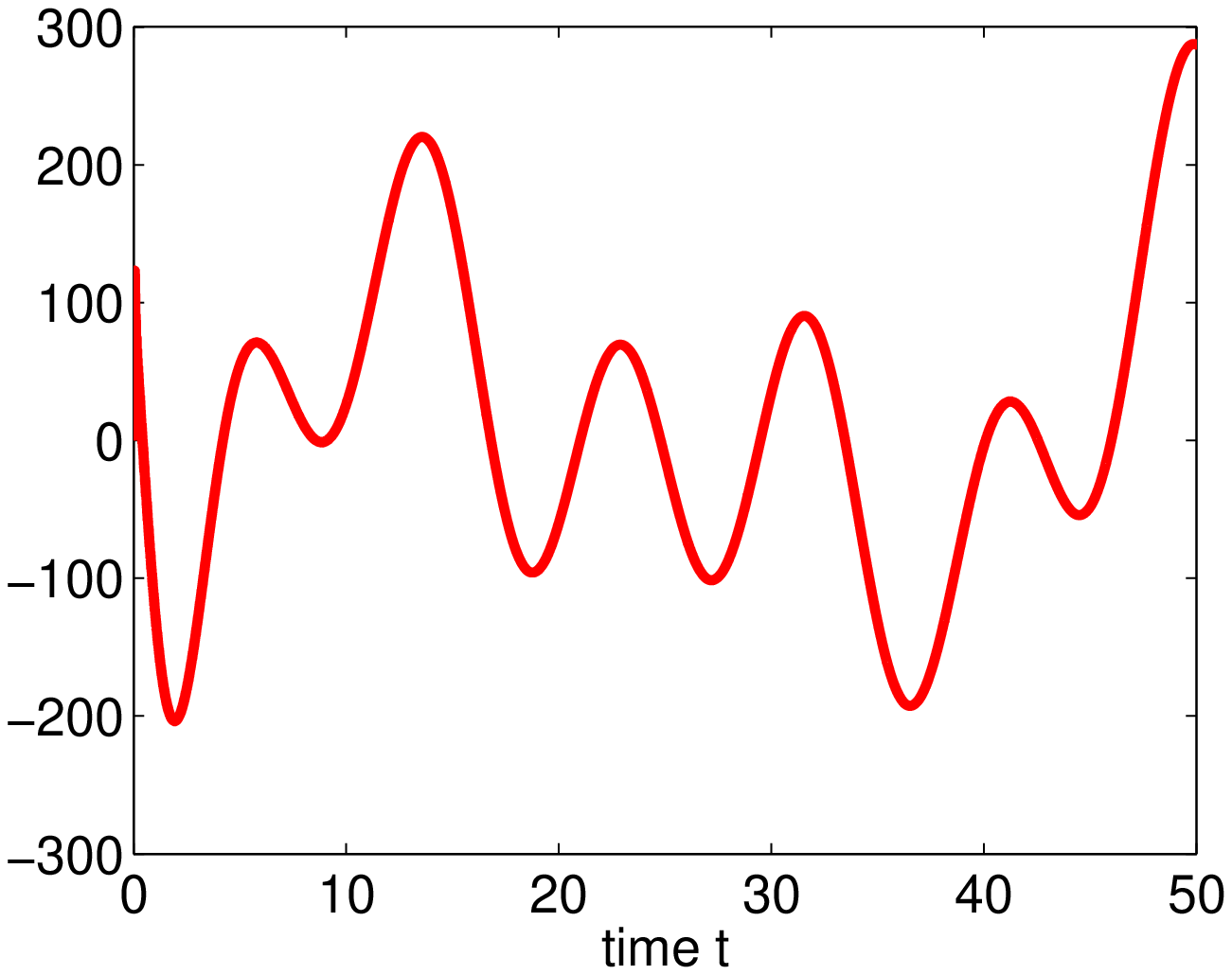,width=2.8in,
height=1.3in}\label{fig:sinu2}}}\caption{Performance of
$\mathcal{L}_1$ adaptive controller with $C(s)=\frac{7500 s
+50^3}{(s+50)^3}$ for $r=100 \cos(0.2 t) $ with time-varying
$\theta(t)=[2+2\cos(0.5 t)\quad 2+0.3\cos(0.5
t)+0.2\cos(t/\pi)]^{\top}$ }
\end{figure}

Next, we consider a higher order filter with low adaptive gain $
\Gamma_c=400\,,\;C(s) = \frac{3 \omega^2 s
+\omega^3}{(s+\omega)^3}\,. $ In Fig. \ref{fig:criterion2}, we
plot
\begin{equation}\label{simu_omega2}
\lambda=\Vert \bar{G}\Vert_{L_1} \theta_{\max}
\end{equation}
 with respect to
$\omega$ and compare it to $1$. We notice that when $\omega
> 25$, we have $\lambda<1$ and  the
$\mathcal{L}_1$-gain requirement in (\ref{condition3}) is
satisfied. Letting $\omega=50$ leads to $ \lambda=0.3984 $.
 The simulation results
of the $\mathcal{L}_1$ adaptive controller are shown in Figs.
\ref{fig:L1_y_r200}-\ref{fig:L1_u_r200} for reference inputs
$r=25,~100,~400$, respectively. We note that it again leads to
scaled control inputs and scaled system outputs for scaled
reference inputs.
In addition, we notice that this performance is achieved by a much
smaller  adaptive gain  as compared to the design with the first
order $C(s)$. Figs. \ref{fig:siny2}-\ref{fig:sinu2} show the
system response and control signal  for reference input $r(t) =
100 \cos (0.2 t)$ and time-varying $\theta(t)=[2+2\cos(0.5 t)\quad
2+0.3\cos(0.5 t)+0.2\cos(t/\pi)]^{\top}$, without any retuning of
the controller.

\begin{rem}
The simulations pointed out that with higher order filter $C(s)$
one could use relatively small adaptive gain. While a rigorous
relationship between the choice of the adaptive gain and the order
of the filter cannot be derived, an insight into this can be
gained from the following analysis. It follows from
(\ref{problem}), (\ref{u1_design}) and (\ref{us_def}) that $ x(s)
= G(s) r(s) +H_o(s) \theta^{\top} x(s) + H_o(s) C(s) \bar{r}(s)\,,
$ while the state predictor in (\ref{G1G2}) can be rewritten as $
\hat{x}(s) = G(s) r(s) + H_o(s) \left( C(s)-1 \right)  \bar{r}(s).
$ We note that $\bar{r}(t)$ is divided into two parts. Its
low-frequency component  $C(s)\bar{r}(s)$ is what the system gets,
while the complementary high-frequency component $\left( C(s)-1
\right) \bar{r}(s)$ goes into the state predictor.
 If the bandwidth of $C(s)$
 is large, then it can suppress only the high frequencies in
 $\bar{r}(t)$, which appear only in the presence of large adaptive
 gain.
 A properly designed higher order $C(s)$
 can be more effective to serve the purpose of filtering with
 reduced
 tailing effects, and, hence
 can generate similar $\lambda$ with smaller bandwidth. This further
 implies that similar performance can be achieved with
  smaller adaptive gain.
\end{rem}

The $\mathcal{L}_1$ adaptive controller has been successfully
flight tested on a miniature aerial vehicle (UAV) with limited
payload capabilities, which consequently restricted the increase
of the adaptation rate \cite{gnc06_mav}. Nevertheless, the flight
tests verified that $\mathcal L_1$ adaptive controller did not
require any tuning. Refs. \cite{gnc06_wingrock,gnc06_airrefueling}
report application of $\mathcal L_1$ controller to different
aerospace benchmark problems.

\section{Conclusion}\label{sec:con}
A novel adaptive control architecture is presented that has
guaranteed transient response in addition to stable tracking. The
new low-pass control architecture adapts fast without generating
high-frequency oscillations in the control signal and leads to
scaled response for both system's input and output signals, which
is otherwise not possible to achieve using conventional adaptive
controllers. The low-frequency behavior of the control signal
implies that the rate saturation is no more an issue. These
arguments enable development of theoretically justified tools for
verification and validation of adaptive controllers. Extension of
the methodology to systems with unknown high frequency gain will
be reported in an upcoming publication.

\end{document}